\documentclass[a4paper,10pt]{amsart}

\usepackage{a4wide,color,eucal,enumerate,mathrsfs}
\numberwithin{equation}{section}

\title[The sticky particle system]
{A Wasserstein approach to the one-dimensional sticky particle system}

\author{Luca Natile}
\address{Dipartimento di Matematica, Universit\`a di Pavia. Via Ferrata, 1 -- 27100 Pavia, Italy.}
\email{luca.natile@unipv.it}

\author{Giuseppe Savar\'e}
\address{Dipartimento di Matematica, Universit\`a di Pavia. Via Ferrata, 1 -- 27100 Pavia, Italy.}
\email{giuseppe.savare@unipv.it}
\urladdr{http://www.imati.cnr.it/\textasciitilde savare}

\thanks{G.S. has been partially supported by MIUR-PRIN'06 grant for the project
 ``Variational methods in optimal mass transportation and in geometric
 measure theory''.}

\keywords{Pressureless Euler equation, Sticky particles, Wasserstein distance, Monotone rearrangement, Gradient flows}

\usepackage{amsthm}
\usepackage{graphicx}
\usepackage{bbm}
\newcommand{\N}{\mathbb{N}}

\newcommand{\R}{\mathbb{R}}

\newcommand{\DD}{\mathscr{D}}

\newcommand{\GG}{\mathscr{G}}

\renewcommand{\SS}{\mathscr{S}}

\newcommand{\calC}{{\ensuremath{\mathcal C}}}

\newcommand{\cH}{{\ensuremath{\mathcal H}}}
\newcommand{\cK}{{\ensuremath{\mathcal K}}}

\newcommand{\cM}{{\ensuremath{\mathcal M}}}
\newcommand{\cN}{{\ensuremath{\mathcal N}}}
\newcommand{\cO}{{\ensuremath{\mathcal O}}}
\newcommand{\cP}{{\ensuremath{\mathcal P}}}
\newcommand{\cT}{{\ensuremath{\mathcal T}}}
\newcommand{\cX}{{\ensuremath{\mathcal X}}}

\newcommand{\cV}{{\ensuremath{\mathcal V}}}

\newcommand{\rrho}{{\mbox{\boldmath$\rho$}}}

\newcommand{\sfi}{{\sf i}}

\newcommand{\sft}{{\sf t}}

\newcommand{\sfv}{{\sf v}}
\newcommand{\sfx}{{\sf x}}

\newcommand{\sfP}{{\sf P}}

\newcommand{\sfS}{{\sf S}}
\newcommand{\sfT}{{\sf T}}

\newcommand{\Kliminf}{K\kern-3pt-\kern-2pt\mathop{\rm lim\,inf}\limits}  
\renewcommand{\d}{{\mathrm d}}
\newcommand{\dt}{{\d t}}

\newcommand{\restr}[1]{\lower3pt\hbox{$|_{#1}$}}
\newcommand{\topref}[2]{\stackrel{\eqref{#1}}#2}
\newcommand{\Leb}[1]{{\mathscr L}^{#1}}      

\newcommand{\up}{\uparrow}
\newcommand{\eps}{\varepsilon}  
\newcommand{\nchi}{{\raise.3ex\hbox{$\chi$}}}
\newcommand{\weakto}{\rightharpoonup}

\newcommand{\media}{\mkern12mu\hbox{\vrule height4pt           %
          depth-3.2pt                                 
          width5pt}\mkern-16.5mu\int\nolimits}        

\newcommand{\Dist}{D}
\newcommand{\Semi}{U}
\newcommand{\Char}{\mathbbm{1}}
\newcommand{\X}{X}
\newcommand{\K}{\cK}
\newcommand{\x}{\tilde x}
\renewcommand{\v}{\tilde v}
\newcommand{\Scalar}[2]{(#1|#2)}
\newcommand{\Rightdert}{{\frac{\d}{\d t}\kern-5pt}^+}
\newcommand{\PRightdert}{{\frac{\partial}{\partial t}\kern-5pt}^+}
\newcommand{\Rightderw}{{\frac{\d}{\d w}\kern-5pt}^+}
\newcommand{\Leftderw}{{\frac{\d}{\d w}\kern-5pt}^-}
\newcommand{\Proj}[1]{\sfP\kern-1pt_{#1}}
\newcommand{\ProjH}[1]{\Proj{\cH_{#1}}}
\newcommand{\band}{\& \relax}
\newcommand{\PositiveCone}{\cN}
\newcommand{\PC}[1]{\PositiveCone_{#1}}

\newcommand{\Xspace}[1]{\cV_{#1}}
\newcommand{\DXspace}{\hat\cV}
\newcommand{\Lspace}[1]{\cX_{#1}(0,1)}
\newcommand{\pseudonorm}[2]{[#2]_{#1}}
\newcommand{\pseudonormp}[2]{\pseudonorm{#1}{#2}^{#1}}
\renewcommand{\sft}{\sfx}

\newtheorem{theorem}{Theorem}[section]

\newtheorem{corollary}[theorem]{Corollary}
\newtheorem{lemma}[theorem]{Lemma}
\newtheorem{proposition}[theorem]{Proposition}

\theoremstyle{definition}
\newtheorem{definition}[theorem]{Definition}
\newtheorem{notation}[theorem]{Notation}

\newtheorem{remark}[theorem]{Remark}

\begin{document}
\begin{abstract}
  We present a simple approach to study the one--dimensional pressureless Euler system
  via adhesion dynamics in the Wasserstein space $\cP_2(\R)$
  of probability measures with finite quadratic moments.

  Starting from a discrete system of a finite number of ``sticky'' particles,
  we obtain new explicit estimates of the solution in terms of the initial
  mass and momentum and we are able to construct an evolution semigroup in
  a measure-theoretic phase space, allowing mass distributions in $\cP_2(\R)$ and
  corresponding $L^2$-velocity fields. We investigate various interesting properties of this semigroup,
  in particular its link with the gradient flow of the (opposite) squared Wasserstein distance.
  
  Our arguments rely on an equivalent formulation of the evolution
  as a gradient flow in the convex cone of nondecreasing functions in the Hilbert space
  $L^2(0,1)$, which corresponds to the Lagrangian system of coordinates given 
  by the canonical monotone rearrangement of the measures.
\end{abstract}
\maketitle
\section{Introduction}
In the recent years considerable attention has been devoted to the
$1$-dimensional pressureless Euler system
\begin{equation}
  \label{eq:Euler}
  \left\{
    \begin{aligned}
      \partial_t \rho + \partial_x (\rho\, v )&=0,\\
      \partial_t (\rho\, v) + \partial_x (\rho\,  v^2)&=0,
    \end{aligned}
    \qquad \text{in }\R\times (0,+\infty);\quad
    \rho\restr{t=0}=\rho_0,\quad
    v\restr{t=0}=v_0,
  \right.
\end{equation}
in connection with the Zeldovich model \cite{Zeldovich70} for the evolution of a ``sticky particle system''
(SPS, in the following) via adhesion dynamics. This model describes the behaviour of
a finite collection of particles, freely moving in absence of
forces and sticking under collision; they can be mathematically represented 
by a time-dependent discrete measure $\rho^N_t:=\sum_{i=1}^n m_i\delta_{x_i(t)}$
concentrated in a finite set of
\begin{quote}
  $N$ particles $P_i(t):=(m_i,x_i(t),v_i(t))$,
  $i=1,\ldots, N$, with positive mass $m_i$, ordered positions
  $x_1(t)\le x_2(t)\le \ldots\le x_{N-1}(t)\le x_N(t)$, and velocities
  $v_i(t)$.
\end{quote}
  Denoting by $J_i(t):=\{j: x_j(t)=x_i(t)\}$ the collection of (the
  indexes of) the particles $P_j(t)$ coinciding with $P_i(t)$ at the
  time $t$, the adhesion dynamic imposes that the sets $J_i(t)$ are
  nondecreasing in time, so that $v_j(t_+)=v_i(t_+)$ for every $j\in
  J_i(t)$.  We can thus order in a finite and monotone sequence
  $0<t_1<t_2<\ldots$ the collection of times when the cardinality of
  some $J_i(t)$ has a discontinuity (corresponding to some collision).
  In each open interval $[t_k,t_{k+1})$ the (right-continuous)
  velocities $v_i(t)=\dot x_i(t)$ are thus supposed to be constant, and
  at each collision time $t_k$ the conservation of mass and momentum
  yields the update equation for the velocities
  \begin{equation}
    \label{eq:1}
    v_i(t_k +) = \frac{\sum_{j\in J_i(t_k)} m_j v_j(t_k-)}
    {\sum_{j\in J_i(t_k         )} m_j}
    ,\qquad 
    i=1,\ldots N.
  \end{equation}
It is not difficult to check that the measures $\rho^N$ and 
$(\rho v)^N_t:=\sum_{i=1}^N m_iv_i(t)\,\delta_{x_i(t)}$ solve \eqref{eq:Euler}.
Starting from the discrete SPS, existence of measure valued solutions to \eqref{eq:Euler} with general initial data and
satisfying suitable entropy conditions
\cite[\textsc{Bouchut}]{Bouchut94}
has been proved by
\textsc{Grenier} \cite{Grenier95} and \textsc{E, Rykov \band Sinai} \cite{E-Rykov-Sinai96}
(but see also the contribution of  \textsc{Martin \& Piasecki} \cite{Martin-Piasecki94})
as limits (in the sense of weak convergence of measures)
of the discrete particle evolutions $\rho^N_t$ as $N\uparrow+\infty$.
Here we also quote 
the different approaches of \textsc{Bouchut \band James} \cite{Bouchut-James95}, of
\textsc{Poupaud \band Rascle} \cite{Poupaud-Rascle97},
and of \textsc{Sever} \cite{Sever01} in the multidimensional case;
viscous regularizations of \eqref{eq:Euler} have been studied by
\textsc{Sobolevski{\u\i}} \cite{Sobolevskii97} and \textsc{Boudin} \cite{Boudin00}, and
a different model, starting from particles of finite size, has been considered by \textsc{Wolansky} \cite{Wolansky07}.

The convergence result has further been extended and refined by \textsc{Brenier \band Grenier} \cite{Brenier-Grenier98},
\textsc{Huang \band Wang} \cite{Huang-Wang01}, and \textsc{Nguyen \band Tudorascu} \cite{Nguyen-Tudorascu08}
(by a different probabilistic approach \textsc{Moutsinga} \cite{Moutsinga08} has recently been able to
consider initial velocities with nonpositive jumps at each points of the support of $\rho_0$):
the basic assumption is 
that the discrete initial velocity $v_i$
is the value in $x_i$ of a given continuous function $v$ with at most linear growth,
and (the total mass being normalized to $1$) the sequence $\rho^N_0$ converges to $\rho_0$
w.r.t.\ the $L^2$-Wasserstein distance in the space $\cP_2(\R)$
of probability measures with finite quadratic moment. This includes the case
(considered in \cite{Brenier-Grenier98})
of a sequence $\rho_0^N$ with uniformly bounded support and weakly converging to $\rho_0$ in the duality
with continuous real functions.

All these results depend on a remarkable characterization of the solution $\rho$ found by
\textsc{Brenier \band Grenier} \cite{Brenier-Grenier98}:
by introducing the cumulative distribution function $M_\rho$ associated to a probability measure $\rho\in \cP(\R)$
\begin{equation}
  \label{eq:53}
  M_{\rho}(x):=\rho((-\infty,x])\quad\forall x\in \R,\quad\text{so that }
  \rho=\partial_x M_{\rho}\quad\text{in }\DD'(\R),
\end{equation}
they prove that the function $M(t,\cdot):=M_{\rho_t}(\cdot)$ is the unique entropy solution of the scalar conservation law
\begin{equation}
  \label{eq:54}
  \partial_t M+\partial_t A(M)=0\quad \text{in }\R\times (0,+\infty),
\end{equation}
where $A:[0,1]\to \R$ is a continuous flux function depending only on $\rho_0$ and $v_0$
(see Theorem \ref{thm:Brenier-Grenier} for a precise statement).

It can also be shown \cite{Nguyen-Tudorascu08} that this solution satisfies the Oleinik entropy condition
\begin{equation}
  \label{eq:55}
  v_t(x_2)-v_t(x_1)\le \frac 1t(x_2-x_1)\quad\text{for $\rho_t$-a.e.\ }x_1,x_2\in \R,\ x_1\le x_2.
\end{equation}
In the present paper we discuss various refinement of \textsc{Brenier-Grenier}
result by a different approach.
Our starting point (Theorem \ref{thm:main1}) is an explicit Lipschitz estimate
(in the $L^p$-Wasserstein distance $W_p$ for every $p\ge 1$, see \eqref{eq:4}) of the dependence of $\rho_t$
with respect to the initial data $\rho_0,(\rho v)_0$: for $p=2$ it shows that
$(\rho_t^N)_{N\in\N}$ is a Cauchy sequence in $\cP_2(\R)$ and in particular
yields the convergence results of \cite{E-Rykov-Sinai96,Brenier-Grenier98,Nguyen-Tudorascu08}
allowing general initial measures in $\cP_2(\R)$ 
and (possibly discontinuous) velocity field $v_0\in L^2(\rho_0)$.
We also show that a suitable $L^2$-like integral distance between the momentum $\rho v$ of two solutions
can be controlled in terms of the initial data and prove further
precise representation properties of the solution and its velocity field (Theorem \ref{thm:precise_sol}).

This leads to the construction of a semigroup $\SS_t$ associated
to the evolution of SPS, which exhibits interesting links with another semigroup
(recently studied by \textsc{Ambrosio, Gigli \band Savar\'e} \cite{Ambrosio-Gigli-Savare08}),
obtained as the gradient flow in $\cP_2(\R)$ of the (opposite) squared Wasserstein distance
from a fixed reference measure.

This link (which  at a first sight may look unexpected) can be better understood in the
simpler case when the initial velocity field $v$ satisfies a one-sided monotonicity condition
(see section 5.4.2 of \textsc{Villani}'s book \cite{Villani03} for more details):
still considering the simpler discrete case, if 
\begin{equation}
  \label{eq:87}
  -\delta^{-1}:=\min_{x_i\neq x_j}\frac {v(x_i)-v(x_j)}{x_i-x_j}<0,\quad
  v(x_i):=v_i,
\end{equation}
for $t\in [0,\delta)$ the map $\sft_0^t(x):=x+tv(x)$ is nondecreasing on the support of $\rho_0$
(the finite set $\{x_i: i=1,\ldots N\}$),
so that the first collision occurs at $t:=\delta $ and in the interval $[0,\delta )$ one has
the freely moving measures
\begin{equation}
  \label{eq:88}
  \rho_t:=(\sft_0^t)_\#\rho_0=\sum_{i=1}^N m_i \delta_{x_i+tv_i},\quad
  (\rho v)_t=\sum_{i=1}^N m_iv_i\delta_{x+tv_i},\ t\in [0,\delta ),
\end{equation}
solving the pressureless Euler system \eqref{eq:Euler}.
On the other hand, the curve $t\mapsto \rho_t$, $t\in [0,\delta ]$, 
is a constant speed minimal geodesic
in $\cP_2(\R)$ connecting $\rho_0$ with $\eta:=\rho_{\delta }$; as in any Riemannian manifold,
it coincides (up to a suitable rescaling, \cite[Theorem 11.2.10]{Ambrosio-Gigli-Savare08})
with the
gradient flow in $\cP_2(\R)$ of the functional $\phi^{\rho_0}(\rho):=-\frac 12W_2^2(\rho,\rho^0)$.
After the collision at time $t=\delta $ the trajectory of the gradient flow does not coincide
with the free motion \eqref{eq:88} anymore, since its velocity has a jump which
can be described exactly by \eqref{eq:1} \cite[Theorem 10.4.12]{Ambrosio-Gigli-Savare08}.
At a later time, the velocity field induced by the (rescaled) Wasserstein gradient flow can be characterized by the
formula
\begin{equation}
  \label{eq:89}
  v_i(t+)=t^{-1}\Big(x_i(t)-\frac{\sum_{j\in J_i(t)} m_j x_j(t)}
    {\sum_{j\in J_i(t)} m_j}\Big)
    ,\qquad 
    i=1,\ldots N,
\end{equation}
and it is an interesting property, stated in Theorems \ref{thm:WGFlow} and \ref{thm:limit}, that
the two different laws \eqref{eq:1} and \eqref{eq:39} give rise to the same evolution, even for arbitrary initial data.

In order to obtain these results, we adopt the point of view of $1$-dimensional
optimal transportation and we represent each probability measures $\rho\in \cP_2(\R)$ by
their monotone rearrangement $X_\rho$, which is the pseudo-inverse of the distribution function $M_\rho$ of \eqref{eq:53}
(a similar approach, in a probabilistic framework, has been also used by \cite{Moutsinga08}; see also
\cite{Gangbo-Nguyen-Tudorascu08} for other applications)
\begin{equation}
  \label{eq:90}
  \X_\rho (w):= \inf
  \left\{x: M_\rho(x) > w \right\}=\inf \big\{x:
  \rho\big((-\infty,x]\big)>w\big\}\quad w\in (0,1).
\end{equation}
The map $\rho\mapsto X_\rho$ is an isometry between $\cP_2(\R)$ (endowed with the $L^2$-Wasserstein distance)
and the convex cone $\K$ of nondecreasing functions in the Hilbert space $L^2(0,1)$. Through this
isometry, any gradient flow with respect to $W_2$ in $\cP_2(\R)$
can be rephrased as a gradient flow in $\K$ with respect to the $L^2(0,1)$-distance, and
one can use the powerful tools the classical theory of variational evolution inequalities in Hilbert spaces 
(we refer to the book by \textsc{Br\'ezis} \cite{Brezis73}).
It turns out (see Theorem \ref{thm:L2flow})
that in this Lagrangian formulation the solution
$X_{\rho_t}$ admits three simple characterizations, in terms of the $L^2(0,1)$-projection
$\Proj\K$ onto $\K$
\begin{equation}
  \label{eq:91}
  X_{\rho_t}=\Proj\K(X_{\rho_0}+tV_0),\quad
  V_0=v_0\circ X_{\rho_0},
\end{equation}
and of the differential inclusions
\begin{equation}
  \label{eq:92}
  \frac \d{\d t}X_{\rho_t}+\partial I_\K(X_{\rho_t})\ni V_0,\qquad
  t\frac \d{\d t}X_{\rho_t}+\partial I_\K(X_{\rho_t})\ni X_{\rho_t}-X_{\rho_0},
\end{equation}
$I_\K$ being the indicator function associated to $\K$ (see next \eqref{eq:22}).
\eqref{eq:91} and \eqref{eq:92} encode all the qualitative information on the measure-valued solution $\rho_t$,
and their proof in the case of the discrete SPS constitutes the core of our argument. It relies on an elementary but careful
description of the $L^2(0,1)$-projection operator $\Proj\K$ and on the subdifferential of $I_\K$, which
has been carried out in Section \ref{sec:K}. Once $\rho_t$ has been determined, its velocity $v_t\in L^2_{\rho_t}(\R)$
can be recovered from the right derivative $V(t):=\Rightdert X_{\rho_t}\in L^2(0,1)$:
in fact, as a byproduct of the second differential inclusion of
\eqref{eq:92}, $V(t)$ is a function of $X(t)$ and therefore one obtains
\begin{equation}
  \label{eq:21}
  V(t)=v_t\circ X_{\rho_t}.
\end{equation}
The projection formula \eqref{eq:91} (which has been introduced by \textsc{Shnirelman} \cite{Shnirelman86,Andrievsky-Gurbatov-Sobolevsky07}
in a slightly different form, see Remark \ref{rem:minimal_lagrangian})
lies more or less explicitly at the core of the formulations by \cite{E-Rykov-Sinai96} and \cite{Brenier-Grenier98}.
As it has been nicely explained by \textsc{Andrievsky, Gurbatov \& Sobolevski{\u\i}} \cite{Andrievsky-Gurbatov-Sobolevsky07}
elaborating the contribution of \cite{Shnirelman86}, \eqref{eq:91} is equivalent to
the \emph{Generalized Variational Principle} of \cite{E-Rykov-Sinai96}, which can be expressed through
the convex envelope of the primitive function of the map $X_{\rho_0}+tV_0$: as stated in full generality
by Theorem \ref{proj on K}, this convexification characterizes the $L^2$-projection on $\K$.
On the other hand, a convexification is also involved in the
second Hopf formula for the solutions of the Hamilton-Jacobi equation associated to \eqref{eq:54},
as it has been already observed by \cite[\S 4]{Brenier-Grenier98}: we will detail this point in
Theorem \ref{thm:Brenier-Grenier}.

The link between the formulation based on the scalar conservation law \eqref{eq:54}
and the Hilbertian theory of gradient flows like \eqref{eq:92}
is not at all surprising, after the illuminating paper by \textsc{Brenier} \cite{Brenier09}.
Wasserstein contraction properties of solutions of one--dimensional scalar conservation laws
have also been recently obtained by \textsc{Bolley, Brenier \band Loeper} \cite{Bolley-Brenier-Loeper05}
(see also the further contribution by \textsc{Carrillo, Di Francesco \band Lattanzio} \cite{Carrillo-DiFrancesco-Lattanzio06}).
So it would be possible in principle to approach SPS starting from \eqref{eq:54} and trying to apply the techniques
developed there. Notice however that two solutions originating from different initial distributions
of position and velocity give rise to two scalar conservation laws differing not only by the initial data but also
by the flux functions, so that their comparison does not look immediate.
Moreover, the present self-contained approach is very simple, since it relies on elementary tools of convex analysis and
direct computations on the discrete case; the simultaneous characterization of the evolution by
\eqref{eq:91} and \eqref{eq:92}
provides a more refined description of the solution and, as a byproduct, a new direct proof of \textsc{Brenier \band Grenier}
theorem.

\subsection*{Plan of the paper}
In the next section we recall some basic definition and notation and we state our main results.
Section \ref{sec:K} collects the main properties related to the convex cone $\K$ in $L^2(0,1)$
(projection, polar cone, subdifferential of the indicator function): they provide simple but crucial tools
for the analysis of the discrete SPS presented in Section \ref{sec:dSPS}, which contains
all the basic calculations.
Section {\ref{sec:lagrangian} deals with existence, stability, and uniqueness of the solution in the Lagrangian formulation.
The final steps of the proofs (mainly concerning the various limit processes)
will be detailed in the last Section \ref{sec:limit}, where we also show
a new derivation of \textsc{Brenier \band Grenier} Theorem \cite{Brenier-Grenier98}
from the Lagrangian representation of the SPS.

\section{Main results}
\label{sec:main}
\subsection*{Couplings, Wasserstein distance, and monotone rearrangementa}
For $p\in [1,+\infty)$ let us denote by $\cP_p(\R)$ the space of Borel probability measures $\rho$
with finite $p$-moment $\int_\R |x|^p\,\d\rho(x)<+\infty$. The $L^p$ Kantorovich-Rubinstein-Wasserstein
distance
$W_p(\rho^1,\rho^2)$ between two measures $\rho^1,\rho^2\in \cP_p(\R)$ can be defined in terms of couplings,
i.e.\ probability measures $\rrho\in \cP(\R\times \R)$  such that $\pi^i_\#\rrho=\rho^i$, $i=1,2$,
by the formula
\begin{equation}
  \label{eq:4}
  W_p^p(\rho^1,\rho^2):=\min\Big\{\int_{\R\times\R}|x-y|^p\,\d\rrho(x,y):\rrho\in \cP(\R\times\R),\
  \pi^i_\#\rrho=\rho^i\Big\}.
\end{equation}
Here $\pi^i(x_1,x_2)=x_i$ is the usual projection on the $i$-th coordinate and
for a general Borel map $\sfT :\R^m\to \R^n$ and a Borel measure $\mu\in \cP(\R^m)$ the
\emph{push-forward} $\nu=\sfT_\#\mu$
is the measure defined by $\nu(A)=\mu(\sfT^{-1}(A))$ for every Borel set $A\subset \R^n$.
We will repeatedly use the change-of-variable formula
\begin{equation}
  \label{eq:30}
  \int_{\R^n}\zeta(y)\,\d(\sfT_\# \mu)(y)=\int_{\R^m}\zeta(\sfT(x))\,\d\mu(x)\quad
  \text{for every Borel map }\zeta:\R^n\to [0,+\infty].
\end{equation}
More generally, given a convex, even, and lower semicontinuous function $\psi:\R\to [0,+\infty]$, we can
consider the cost $c_\psi(x,y):=\psi(x-y)$, $x,y\in \R$, and the associated optimal mass transportation problem
\begin{equation}
  \label{eq:97}
  \calC_\psi(\rho^1,\rho^2):=\inf\Big\{\int_{\R\times\R}\psi(x-y)\,\d\rrho(x,y):\rrho\in \cP(\R\times\R),\
  \pi^i_\#\rrho=\rho^i\Big\}.
\end{equation}
In the present $1$-dimensional case, there exists a unique
optimal coupling $\rrho=\Gamma_o(\rho^1,\rho^2)$ realizing the minimum
of \eqref{eq:4} and of \eqref{eq:97} (at least when the cost is finite):
it can be explicitly characterized by inverting the distribution functions of
$\rho^1,\rho^2$. More precisely, for every $\rho\in \cP(\R)$ we consider
its monotone rearrangement $X_\rho$ \eqref{eq:90},
a right-continuous and nondecreasing function satisfying
\begin{equation}
  \label{eq:6}
  (\X_\rho)_\#\lambda=\rho,\quad
  \lambda:=\Leb 1\restr{(0,1)},\qquad
  \int_\R \zeta(x)\,\d\rho(x)=\int_0^1 \zeta(X_\rho(w))\,\d w
\end{equation}
for every nonnegative Borel map $\zeta:\R\to [0,+\infty]$.
In particular,
$\rho\in \cP_p(\R)$ iff $\X_\rho \in L^p(0,1)$. Moreover, thanks to the Hoeffding-Fr\'echet theorem
\cite[Sec. 3.1]{Rachev-Ruschendorf98I},
the joint map $\X_{\rho^1,\rho^2}(w):=(\X_{\rho^1}(w),
\X_{\rho^2}(w)),$ $w\in (0,1)$, characterizes the optimal coupling $\rrho\in \Gamma_o(\rho^1,\rho^2)$
by the formula 
\begin{equation}
  \label{eq:8}
  \rrho=\big(\X_{\rho^1,\rho^2}\big)_\# \lambda,
\end{equation}
so that \cite{DallAglio56,Rachev-Ruschendorf98I,Villani03}
\begin{equation}
  \label{W_2}
  W_p^p(\rho^1, \rho^2) = \int_0^1 \big|\X_{\rho^1}(w)-\X_{\rho^2}(w)\big|^p\,\d w,\quad
  \calC(\rho^1,\rho^2)=\int_0^1 \psi\big(\X_{\rho^1}(w)-\X_{\rho^2}(w)\big)\,\d w,
\end{equation}
and the map $\rho\in \cP(\R)\longmapsto \X_\rho$
is an isometry between $\cP_2(\R)$ and the convex subset $\K$ of $L^2(0,1)$ of (essentially)
nondecreasing functions
(which can be identified with their right-continuous representatives).

\subsection*{An explicit estimate through Wasserstein distance}
We introduce the set
\begin{equation}
  \label{eq:10}
  \Xspace p (\R):=\Big\{\mu=(\rho,\rho v)\in \cP_p(\R)\times \cM(\R):  v\in L^p_\rho(\R)\Big\},\quad p\in [1,+\infty),
\end{equation}
$\cM(\R)$ being the set of all signed Borel measures with finite total variation,
the semi-distances (here $\mu^i=(\rho^i,\rho^i v^i)$)
\begin{align}
  \label{eq:11}
  \Semi_p^p(\mu^1,\mu^2):=&\int_{\R\times\R}|v^1(x)-v^2(y)|^p\,\d\rrho(x,y)\quad
  \rrho=\Gamma_o(\rho^1,\rho^2)\\
  \label{eq:12}
  =&   \int_0^1 |v^1(\X_{\rho^1}(w))-v^2(\X_{\rho^2}(w))|^p\,\d w,
\end{align}
and the distances
\begin{align}
  \label{eq:11bis}
  \Dist^p_p(\mu^1,\mu^2):=&W_p^p(\rho^1,\rho^2)+\Semi_p^p(\mu^1,\mu^2).
\end{align}
We also set
\begin{equation}
  \label{eq:67}
  \pseudonormp p{\mu}:=\int_\R \Big(|x|^p+v^p(x)\Big)\,\d\rho(x)=
  \Dist_p^p\big(\mu,(\delta_0,0)\big).
\end{equation}
\begin{proposition}
  \label{prop:metric}
  $\Dist_p$ is a distance in $\Xspace p(\R)$ and $(\Xspace p (\R),\Dist_p)$ is metric (but not complete)
  space whose topology is stronger
  than the one induced by the weak convergence of measures.
  The collection of discrete measures
  \begin{equation}
    \label{eq:15}
    \DXspace(\R):=\Big\{\mu=\big(\sum_{i=1}^N m_i \delta_{x_i},\sum_{i=1}^N
    m_i v_i\delta_{x_i}\big): m_i>0,\ \sum_{i=1}^N m_i=1,\ x_i,v_i\in \R\Big\}
  \end{equation}
  is a dense subset of $\Xspace p (\R)$.
  A sequence $\mu_n=(\rho_n,\rho_n v_n)$, $n\in \N$, converges to $\mu=(\rho,\rho v)$ in
  $\Xspace p(\R) $, $p>1$, if and only if (see \cite[Def. 5.4.3]{Ambrosio-Gigli-Savare08})
  \begin{equation}
    \label{eq:14}
    W_p(\rho_n,\rho)\to0,\quad \rho_nv_n\weakto \rho v\quad \text{weakly in }\cM(\R),\quad
    \int_\R |v_n|^p\,d\rho_n\to \int_\R |v|^p\,\d\rho.
  \end{equation}
\end{proposition}
Let us denote by $\SS_t:\DXspace(\R)\to \DXspace(\R)$ the map associating to
any discrete initial datum $(\rho_0,\rho_0v_0)$ the solution $(\rho_t,\rho_tv_t)$
of the (discrete) sticky-particle system. $\SS_t$ is a semigroup in $\DXspace(\R)$.
\begin{theorem}[Stability with respect to the initial data]
  \label{thm:main1}
  Let $\mu^\ell_t=(\rho^\ell_t,\rho^\ell_tv^\ell_t)=\SS_t[\mu^\ell_0]$, $\ell=1,2$,
  be the solutions of the (discrete) sticky-particle system
  with initial data $\mu_0^\ell\in \DXspace(\R)$.
  Then for every convex cost \eqref{eq:97} and every $p\ge 1$
  \begin{subequations}
    \label{sub:16}
    \begin{align}
      \label{eq:98}
      \calC_\psi(\rho^1_t,\rho^2_t)&\le \int_{\R\times\R}\psi\big(x+tv^1(x)-(y+tv^2(y)\big)\,\d\rrho(x,y),
      \quad \rrho=\Gamma_o(\rho^1,\rho^2),\\
      \label{eq:16}
      W_p(\rho^1_t,\rho^2_t)&\le W_p(\rho^1_0,\rho^2_0)+t\Semi_p(\mu^1_0,\mu^2_0),\\
      \label{eq:17}
      \int_0^t\Semi^2_2(\mu^1_r,\mu^2_r)\,\d r&\le
      C(1+t)\Big(\pseudonorm2{\mu^1}+\pseudonorm2{\mu^2}\Big)
      \Big(W_2(\rho^1_0,\rho^2_0)+\Semi_2(\mu^1_0,\mu^2_0)\Big),
    \end{align}
  \end{subequations}
  for a suitable ``universal'' constant $C$ independent of $t$ and the data.
 \end{theorem}
 We say that a map $\SS:\Xspace p(\R)\to \Xspace p(\R)$ is \emph{strongly-weakly (s-w) continuous}
 if for every $\mu^n,\mu\in \Xspace p(\R)$ with $\SS[\mu^n]=
 (\tilde \rho^n,\tilde \rho^n\tilde v^n),\ \SS[\mu]=(\tilde\rho,\tilde\rho \tilde v)\in \Xspace p(\R),$
 \begin{equation}
   \label{eq:111}
   \lim_{n\uparrow+\infty}\Dist_p(\mu_n,\mu)=0\quad
   \Longrightarrow\quad
   \lim_{n\uparrow+\infty}W_p(\tilde\rho^n,\tilde\rho)=0,\quad
   \tilde\rho_n\tilde v_n\weakto \tilde\rho \tilde v\quad\text{weakly in }\cM(\R).
 \end{equation}
\begin{theorem}[The evolution semigroup in $\Xspace p(\R)$]
  \label{thm:precise_sol}
  \ 
  \begin{enumerate}[(a)]
  \item 
    The semigroup $\SS_t$ can be uniquely extended by density to a
    right-continuous semigroup (still denoted $\SS_t$) of strongly-weakly continuous
    transformations in $\Xspace p(\R)$, $p\ge 2$, thus satisfying
    \begin{equation}
      \label{eq:110}
      \SS_{s+t}[\mu]=\SS_s[\SS_t[\mu]]\quad\forall\, s,t\ge0,\qquad
      \lim_{t\downarrow0}\Dist_p(\SS_t[\mu],\mu)=0\quad
      \forall\, \mu\in\Xspace p(\R).
    \end{equation}
    $\SS_t$ complies with the same estimates $(\ref{sub:16}a,b,c)$
    of Theorem \ref{thm:main1}.
    \item $(\rho_t,\rho_tv_t)=\SS_t[\mu]$, $\mu\in \Xspace2(\R)$,
      is a distributional solution of \eqref{eq:Euler} satisfying
      Oleinik entropy condition \eqref{eq:55}.
    \item If $\psi:\R\to\R$ is a convex function such that $\psi(v_0)\in L^1_{\rho_0}(\R)$,
      and $(\rho_t,\rho_t v_t)=\SS_t[\mu_0]$, then 
      \begin{equation}
        \label{eq:100}
        \text{the map}\quad
        t\mapsto \int_\R \psi(v_t)\,\d\rho_t(x)
        \quad
        \text{is nonincreasing in $[0,+\infty)$},
      \end{equation}
      and its (at most countable) jump set $\cT=\cT(\mu)$ is independent of $\psi$.
    \item
      If $\mu\in \Xspace p(\R)$ and
      $\mu_t=(\rho_t,\rho_tv_t)=\SS_t[\mu]$, $t\in [0,+\infty)$, the
      curve $t\mapsto \rho_t$ is Lipschitz in $\cP_p(\R)$ with respect
      to $W_p$, and the curve $t\mapsto \rho_tv_t$ is continuous with
      respect to the weak topology in $\cM(\R)$, right-continuous in $[0,+\infty)$ with
      respect to the (semi-) distance $\Semi_p$, and left-continuous at each
      $t\in (0,+\infty)\setminus \cT$ where $\cT$ is the at most
      countable jump set of \eqref{eq:100}.
    \item
      Let $\mu^n_t=(\rho^n_t,\rho^n_t v^n_t)=\SS_t[\mu^n]$ and $\mu_t=(\rho_t,\rho_tv_t)=\SS_t[\mu]$;
      if $\mu^n$ converges to $\mu$ in $\Xspace p(\R)$ as $n\uparrow+\infty$,
      then for every $t\in [0,+\infty)$ $\rho^n_t$ converges to $\rho_t$ in $\cP_p(\R)$,
      $\rho^n_tv^n_t$ weakly converges to $\rho_t v_t$ in $\cM(\R)$;
      moreover, $\mu^n_t$ converges to $(\rho_t,\rho_t v_t)=\SS_t[\mu]$ in $\Xspace p(\R)$ for every
      $t\in [0,+\infty)\setminus \cT(\mu)$.
    \item
      For every $0\le s<t$
      there exists a $\rho_s$-essentially unique monotone map
    $\sft_s^t\in L^2_{\rho_s}(\R)$ such that
    \begin{equation}
      \label{eq:71}
      \rho_t=(\sft_s^t)_\#\rho_s,\quad
      \lim_{h\downarrow 0}\frac{\sft_s^{s+h}-\sfi}h=v_s\quad
      \text{in }L^2_{\rho_s}(\R),\quad \sfi(x)\equiv x,
    \end{equation}
    \begin{equation}
      \label{eq:72}
      v_t(y)=\int_\R v_s(x)\,\d\rho^{s\to t}_y(x)=
      (t-s)^{-1}\Big(y-\int_\R \sft_s^t(x)\,\d\rho^{s\to t}_y(x)\Big)\quad\text{for $\rho_t$-a.e.\ $y\in \R,$}
    \end{equation}
    where $\rho^{s\to t}_y$ is the disintegration of $\rho_s$ with
    respect to $\sft_s^t$.
  \end{enumerate}
\end{theorem}
Let us recall that the disintegration $\rho^{s\to t}_y$ of $\rho_s$ with respect to
the Borel (monotone) map $\sft_s^t$ is a Borel family of
parametrized measures uniquely determined for $\rho_t$-a.e.\ $y\in \R$, such that
$\rho_s=\int_\R \rho^{s\to t}_y\,\d\rho_t(y)$ with $\rho^{s\to t}_y((\sft_s^t)^{-1}(y))=1$
(see e.g. \cite[Thm. 5.3.1]{Ambrosio-Gigli-Savare08}).

Notice that for a fixed $t$ the map $\SS_t:\Xspace p(\R)\to \Xspace p(\R)$ may fail to be continuous
with respect to the distance $\Dist_p$, at least in the momentum component $\rho v$. 

\subsection*{The gradient flow of the (opposite) squared Wasserstein distance}
\eqref{eq:71} and \eqref{eq:72} show an interesting connection between the semigroup $\SS_t$ in $\Xspace2 (\R)$
and the gradient flow $\GG_t^\sigma$ in $\cP_2(\R)$ of the (opposite) squared distance functional
\begin{equation}
  \label{eq:20}
  \phi^{\sigma}(\rho):=-\frac 12 W_2^2(\rho,\sigma)\quad \forall\, \rho,\sigma\in \cP_2(\R).
\end{equation}
Let us recall \cite{Ambrosio-Gigli-Savare08} that for every choice of a reference measure $\sigma\in \cP_2(\R)$
it is possible to define a unique continuous and $1$-expansive semigroup $\GG^{\sigma}_\tau:
\cP_2(\R)\to \cP_2(\R)$, $\tau\ge 0$, whose Lipschitz trajectories $\hat \rho_\tau:=\GG^{\sigma}_\tau(\rho)$
can be uniquely characterized by the Evolution Variational Inequality
\begin{equation}
  \label{eq:24tau}
  \frac 12\frac \d{\d \tau} W_2^2(\hat \rho_\tau,\eta)-
  \frac 12 W_2^2(\hat\rho_\tau,\eta)\le
  \phi^{\sigma}(\eta)- \phi^{\sigma}(\hat \rho_\tau)\quad
  \forall\, \eta\in \cP_2(\R).
\end{equation}
The next result shows that $\SS_t$ and $\GG^{\rho_0}_\tau$ basically coincide, up to the rescaling
\begin{equation}
  \label{eq:73}
  \tau=\log t,\quad t=e^\tau,\quad \hat\rho_\tau=\rho_{e^\tau}.
\end{equation}
\begin{theorem}[Gradient flow of the Wasserstein distance and SPS]
  \label{thm:WGFlow}
  Let $(\rho_t,\rho_tv_t)=\SS_t(\rho_0,\rho_0v_0)\in \Xspace2 (\R)$ be the semigroup
  solution of the sticky-particle system.
  The Lipschitz curve $(\rho_t)_{t\ge 0}$ in $\cP_2(\R)$
  for a.e.\ $t>0$ it solves
  the Evolution Variational Inequality
  \begin{equation}
    \label{eq:24}
    \frac t2\frac \d{\d t} W_2^2(\rho_t,\eta)-\frac 12 W_2^2(\rho_t,\eta)\le
    \phi^{\rho_0}(\eta)-\phi^{\rho_0}(\rho_t)\quad\text{a.e.\ in }(0,+\infty),\quad
    \forall\, \eta\in \cP_2(\R).
  \end{equation}
  Equivalently, the reparametrized solutions $\hat\rho_\tau=\rho_{e^\tau}$ satisfy
  \eqref{eq:24tau} with $\sigma:=\rho_0$ and 
  we thus get the representation formula
\begin{equation}
  \label{eq:28}
  \hat\rho_\tau=\GG^{\rho_0}_{\tau-\delta}\hat \rho_\delta\quad
  \text{or, equivalently,}\quad
  \rho_t=\GG^{\rho_0}_{\log (t/\eps)}\rho_\eps\quad\forall\, \tau=\log t\ge\delta=\log \eps.
\end{equation}
Conversely, if $t\mapsto \rho_t$ is a Lipschitz curve in $\cP_2(\R)$ satisfying \eqref{eq:24} and
the initial velocity condition
\begin{equation}
  \label{eq:109}
  \lim_{t\downarrow0}t^{-2}\int_\R |x+tv_0(x)-y|^2\,\d\rrho_t(x,y)=0\quad \rrho_t=\Gamma_o(\rho_0,\rho_t),
\end{equation}
then there exists a unique Borel velocity vector field $v_t\in L^2_{\rho_t}(\R)$ such that
$(\rho_t,\rho_tv_t)=\SS_t(\rho_0,\rho_0 v_0)$. $v_t$ is the Wasserstein velocity field of $\rho_t$
\cite[Thm. 8.4.5]{Ambrosio-Gigli-Savare08}.
\end{theorem}
Notice that \eqref{eq:109} corresponds to \eqref{eq:71} for $s=0$
in the case (which a posteriori is always verified) $\rrho_t=(\sfi\times \sft_0^t)_\#\rho_0$.

We can use \eqref{eq:28} to exhibit the solution $\rho_t$ of SPS by a
simple limit procedure:
\begin{theorem}
  \label{thm:limit}
  Let $(\rho_t,\rho_t v_t)=\SS_t(\rho_0,\rho_0v_0)\in \Xspace2 (\R)$ be the solution of SPS
  and let $\tilde\rho_\eps:=(i+\eps v_0)_\# \rho_0$, $\eps>0$.
  Then
  \begin{equation}
    \label{eq:29}
    \rho_t=\lim_{\eps\downarrow0} \GG^{\rho_0}_{\log(t/\eps)}(\tilde\rho_\eps)\quad\text{in }\cP_2(\R).
  \end{equation}
  Moreover, if for some $\eps_0>0$ the map $i+\eps_0 v_0$ is $\rho_0$-essentially nondecreasing then
  \begin{equation}
    \label{eq:19}
    \rho_\eps=\tilde\rho_\eps,\quad
    \rho_t= \GG^{\rho_0}_{\log(t/\eps)}(\tilde\rho_\eps)\quad \forall\, \eps\in (0,\eps_0],\ t\ge \eps.
  \end{equation}
\end{theorem}
\subsection*{The evolution in Lagrangian coordinates}
We conclude this section with an even more explicit formula for the evolution of
the monotone rearrangement function $\X(t)=\X_{\rho_t}$.
We denote by $I_\K$ the indicator (convex, lower semicontinuous) function of $\K$ in $L^2(0,1)$
\begin{equation}
  \label{eq:101}
  I_{\K}(\X)=
  \begin{cases}
     0&\text{if }\X\in \K,\\
    +\infty&\text{otherwise,}
  \end{cases}
  \quad \text{with subdifferential}\quad
  \partial I_\K:L^2(0,1)\to 2^{L^2(0,1)}.
\end{equation}
We also introduce the closed subspace $\cH_X\subset L^2(0,1)$, $X\in \K$, whose functions $Y\in L^2(0,1)$
are essentially constant in each open interval $(a,b)\subset (0,1)$ where $X$ is constant:
it is not difficult to check that for every $X\in \K$ and $Y\in L^2(0,1)$
\begin{equation}
  \label{eq:116}
  \text{$Y\in \cH_X$\quad iff\quad
    $Y=y\circ X$ for some Borel map $y\in L^2_\rho(\R)$, $\rho=X_\#\lambda$.}
\end{equation}
\begin{theorem}[Lagrangian evolution]
  \label{thm:L2flow}
  A curve $(\rho_t,\rho_tv_t)\in \Xspace2(\R)$, $t\ge0$, is the semigroup solution $\SS_t(\rho_0,\rho_0v_0)$ of
  SPS
  as in Theorem \ref{thm:main1} if and only if
  its monotone rearrangement $X(t)=X_{\rho_t}\in \K\subset L^2(0,1)$ satisfies one of the following three (equivalent)
  characterizations in terms of the couple 
  $X_0:=X_{\rho_0}$ and $V_0:=v_0(X_0)\in \cH_{X_0}$:
  \begin{enumerate}[{\bf I.}]
  \item  $X$ is the unique strong (i.e.\ absolutely continuous) solution  of the Cauchy
  problem for the subdifferential inclusion
  \begin{equation} \label{eq:44} \frac{\d}{\d t} \X \in - \partial
    I_\K(\X)+V_0,\quad X(0)=X_0.
    \tag{L.I}
  \end{equation}
  \item  $X$ admits the representation formula
  \begin{equation}
    \label{eq:49}
    \X(t)= \Proj \K(\X_0 + t V_0)
    \tag{L.II}
  \end{equation}
  where $\Proj \K$ is the $L^2$-projection on the convex cone
  $\K\subset L^2(0,1)$.
  \item
    $X$ is the unique strong solution of the rescaled
    gradient flow
    \begin{equation}
      \label{eq:26}
      t\frac \d{\d t}\X(t)\in -\partial I_{\K}(\X(t))+\X(t)-\X_0,\quad\text{such that}\quad
      \lim_{t\downarrow 0}t^{-1}(X(t)-X_0)=V_0\quad\text{in }L^2(0,1).
      \tag{L.III}
    \end{equation}
  \end{enumerate}
  In each of these cases the curve $t\mapsto \X(t)$ is Lipschitz continuous in $L^2(0,1)$ and right-differentiable
  at each time $t$; the velocity field $v_t$ can be recovered by the formula 
  \begin{equation}
    \label{eq:74}
    V(t)=\Rightdert X(t)=v_t\circ X(t)=\Proj{\cH_{X(t)}}(V_0)\in \cH_{X(t)}
    \quad\forall\, t\ge 0,
    \tag{L.a}
  \end{equation}
  where $\Proj {\cH_X}$ denotes the $L^2$ orthogonal projection on the closed subspace
  $\cH_X\subset L^2(0,1)$.
  The closed subspaces $\cH_{X(t)}$ are nonincreasing
  \begin{equation}
    \label{eq:108}
    \cH_{X(t)}\subset \cH_{X(s)}\quad \text{if }0\le s\le t,
    \tag{L.b}
  \end{equation}
  and $X,V$ satisfy the semigroup identities
  \begin{align}
    \label{eq:49bis}
    \X(t)= \Proj\K(X(s)+(t-s) V(s)),\quad
    V(t)=\Proj{\cH_{X(t)}}(V(s))\quad
    \forall\, 0\le s\le t.
    \tag{L.c}
  \end{align}
\end{theorem}
This result shows that the natural evolution space
for the Lagrangian sticky particles flow is
\begin{equation}
  \label{eq:113}
  \Lspace p:=\Big\{(X,V)\in L^p(0,1)\times L^p(0,1): X\in \K,\  V=v\circ X\in \cH_X\Big\}\quad p\ge 2,
\end{equation}
endowed with the product distance in $L^p(0,1)\times L^p(0,1)$. The bijective map
\begin{equation}
  \label{eq:5}
  (\rho,\rho v)\in \Xspace p(\R)\longleftrightarrow (X,V)\in \Lspace p,\quad X=X_\rho,\ V=v\circ X_\rho,
\end{equation}
is in fact an isometry with respect to $\Dist_p$ of \eqref{eq:11bis}.
\begin{corollary}[Lagrangian semigroup]
  \label{cor:semigroup}
  For every $p\ge 2$
  the time dependent transformations $\sfS_t:\Lspace p\to \Lspace p$, $t\ge0$, which map
  a couple $(X_0,V_0)\in \Lspace p$ into the couple $(X(t),V(t))=\sfS_t(X_0,V_0)\in \Lspace p$
  where $X$ is the solution of (one of the equivalent) \emph{(L.I, II, III)} and
  $V=\Rightdert X$ as in \eqref{eq:74}, define a right-continuous semigroup in $\Lspace p$,
  satisfying
  \begin{equation}
    \label{eq:114}
    \begin{gathered}
      (X(t),V(t))=\sfS_t(X_0,V_0)\quad \Longleftrightarrow\quad
      (\rho_t,\rho_t v_t)=\SS_t(\rho_0,\rho_0v_0)\\
      \text{where}\qquad
      \rho_t=\big(X(t)\big)_\#\lambda,\quad V(t)=v_t\circ X(t).
    \end{gathered}
  \end{equation}
\end{corollary}
\begin{remark}[Rescaling]
  \label{rem:rescaling}
  Up to the rescaling $\tau=\log t$, $\hat X(\tau)=X(e^\tau)$,
  \eqref{eq:26} is equivalent to
  \begin{equation}
    \label{eq:26tau}
    \frac \d{\d \tau }\hat \X(t)\in 
    -\partial I_{\K}(\hat \X(t))+\hat \X(t)-\X_0.
  \end{equation}
\end{remark}
We shall show (see Theorem \ref{proj on K}) the $\Proj\K$ is a contraction in every $L^p(0,1)$, so that
\eqref{eq:49} provides a simple and sharp way to estimate $X(t)$ in terms of the initial data
corresponding to \eqref{eq:16}. Applying a general result of
\cite{Savare93,Savare96}, one can obtain \eqref{eq:17} from the representation
\eqref{eq:44}.

Let us finally remark that the Wasserstein gradient flow of Theorem \ref{thm:WGFlow}
is equivalent to \eqref{eq:26}-\eqref{eq:26tau}: it is sufficient to introduce the functional $\Phi^\sigma$
\begin{equation}
  \label{eq:22}
  \Phi^\sigma(\X):= -\frac 12 \|\X-\X_\sigma\|_{L^2(0,1)}+I_{\K}(\X),\quad
  X\in L^2(0,1),
\end{equation}
which is related to $\phi^\sigma$ by
\begin{equation}
  \label{eq:23}
  \phi^\sigma(\rho)=\Phi^\sigma(\X_\rho)\quad
  \forall\, \rho\in \cP_2(\R),
\end{equation}
and is a smooth quadratic perturbation of the convex and lower semicontinuous
indicator functional $I_{\K}$; since 
\begin{equation}
  \label{eq:25}
  \partial \Phi^\sigma(\X)=\partial I_{\K}(\X)-(\X-\X_\sigma),
\end{equation}
\eqref{eq:26tau} is the subdifferential formulation in $L^2(0,1)$ of the gradient flow of $\Phi^{\rho_0}$,
whose metric characterization \cite{Ambrosio-Gigli-Savare08} yields \eqref{eq:24tau} thanks to the isometry
$\rho\leftrightarrow X_\rho$
between $\cP_2(\R)$ and $\K$.
\begin{remark}[Minimal Lagrangian description]
  \label{rem:minimal_lagrangian}
  One can use (as in \cite{Shnirelman86,Andrievsky-Gurbatov-Sobolevsky07})
  the initial measure $\rho_0\in \cP(\R)$ as
  a reference for the Lagrangian evolution, thus representing $\rho_t$ as $\sft(t)_\#\rho_0$
  for the optimal monotone map $\sft(t)=\sft_0^t
  \in L^2_{\rho_0}(\R)$ according to Theorem \ref{thm:precise_sol} (f). 
 We can therefore introduce the convex set $\K(\rho_0)$ of essentially nonincreasing
  Borel maps in the Hilbert space $L^2_{\rho_0}(\R)$ and we have the corresponding formulae
  for the evolution in $L^2_{\rho_0}(\R)$ ($\sfi:\R\to\R$ denotes the identity map)
  \begin{gather}
    \label{eq:112}
    \frac\d{\dt}\sft(t)\in -\partial I_{\K(\rho_0)}(\sft(t))+v_0,\quad \sft(0)=\sfi,
    \tag{L.I'}\\
    \sft(t)=\Proj{\K({\rho_0})}(\sfi+tv_0),\quad \sfi(x)=x,
    \tag{L.II'}\\
    t\frac \d{\dt}\sft(t)\in - \partial I_{\K({\rho_0})}(\sft(t))+\sft(t)-\sfi,
    \tag{L.III'}
  \end{gather}
  to be completed with the expression for the velocity $\Rightdert \sft(t)=\sfv(t)=v_t\circ\sft(t)$.
  All these relations could be easily deduced by
  Theorem \ref{thm:L2flow}, since the correspondence
  $\sfx\leftrightarrow X=\sfx\circ X_0$ is an isometry between $L^2_{\rho_0}(\R)$ and the closed subspace
  $\cH_{X_0}$ of $L^2(0,1)$. On the other hand, it is easier to deal with the convex set $\K$
  in the space $L^2(0,1)$ with the uniform Lebesgue measure as a reference and the description
  provided by Theorem \ref{thm:L2flow} is more general, since it allows to compare
  solutions arising from different initial data.
\end{remark}
\section{Main properties of $\K$}
\label{sec:K}
In this section we will study the properties of the convex set $\K$
of nondecreasing functions in $L^2(0,1)$, in particular the $L^2(0,1)$-projection operator $\Proj \K$
and the subdifferential of the indicator function $I_\K$ \eqref{eq:101}.
Denoting by $\Scalar\cdot\cdot$ (resp. $\|\cdot\|$) the usual scalar product (resp. the induced norm) in $L^2(0,1)$,
since $\K$ is a convex cone, $\Proj \K$ can be characterized by
\begin{align}
  \label{eq:34}
  g=\Proj \K(f)\quad &\Longleftrightarrow\quad
  g\in \K,&
  \Scalar{f-g}{z-g}\le 0\quad &\forall\, z\in \K
   \\\label{eq:34bis}
   &\Longleftrightarrow\quad
   g\in \K,&
   \quad \Scalar{f-g}z \le 0\quad  &\forall\, z \in \K,
   \qquad
   \Scalar{f-g}g= 0.  
\end{align}
The next result provides a useful characterization of $\Proj \K(f)$
in terms of the convex envelope of the primitive of $f$.
Recall that the convex envelope of a given continuous 
function $F:[0,1]\to \R$ is defined as
\begin{equation}
  \label{eq:31}
  F^{**}(w):=\sup\Big\{a+bw:a,b\in \R,\ a+b v\le F(v)\quad \forall\, v\in [0,1]\Big\}\quad w\in [0,1],
\end{equation}
and it is the greatest bounded, (lower semi-) continuous, and convex function $G$ satisfying $G\le  F$ in $[0,1]$; it is therefore
right and left differentiable at every point $t\in (0,1)$ and its right derivative $g:=\Rightderw F^{**}$
is nondecreasing and right continuous.
\begin{theorem}[Projection on $\K$]
  \label{proj on K}
  Let $f \in L^2(0,1)$ and let $F(w) = \int_0^w f(s) ds$ be its primitive. Then
  $$
  \Proj \K(f) =g=\Rightderw F^{**}
  $$
  where $F^{**}$ is the convex envelope of $F$ defined by \eqref{eq:31}.
  Moreover, for every convex lower semicontinuous function $\psi:\R\to (-\infty,+\infty]$
  and every $f,h\in L^2(0,1)$ we have
  \begin{equation}
    \label{eq:96}
    \int_\R \psi\big(\Proj\K(f)\big)\,\d w\le \int_\R \psi(f)\,\d w,\quad
    \int_\R \psi\big(\Proj\K(f)-\Proj\K(h)\big)\,\d w\le \int_\R \psi(f-h)\,\d w.
  \end{equation}
  In particular, $\Proj\K$ is a contraction in every space $L^p(0,1)$, $p\in [1,+\infty]:$
  \begin{equation}
    \label{eq:102}
    \|\Proj\K(f)-\Proj\K(h)\|_{L^p(0,1)}\le \|f-g\|_{L^p(0,1)}\quad \forall\, f,h\in L^p(0,1).
  \end{equation}
\end{theorem}
We split the proof in several steps.
Here is a preliminary Lemma.
\begin{lemma} \label{le:derivatives}
For every $f\in L^2(0,1)$ $F^{**}$ is continuous in $[0,1]$, locally Lipschitz in $(0,1),$ and coincides with $F$ 
at $w=0$ and $w=1$. If $f\in L^\infty(0,1)$ then $F$ and $F^{**}$ are Lipschitz continuous in the closed interval $[0,1]$.
\end{lemma}
\begin{proof}
  Let us first assume $f\in L^\infty(0,1)$ and let $L$ be the Lipschitz constant of $F$; then
  $$
  F(0) - L w \le F(w),\quad F(1) + L(w-1) \le  F(w)\quad \forall\, w\in [0,1],
  $$
  so that $F^{**}(0)=F(0)$, $F^{**}(1)=F(1)$, and 
  \begin{equation} \label{ciao}
    F(0) - Lw \le F^{**}(w),\quad  F(1) + L(w-1) \le F^{**}(w)\quad \forall\, w\in [0,1].
  \end{equation}
  Therefore the right derivative $g$ of $F^{**}$ satisfies 
  $-L \le g(0) \le g(w)\le \Leftderw F^{**}(1)\le L$ so that $F^{**}$ is a Lipschitz function.

  In the general case when $f\in L^2(0,1)$, we can approximate its (absolutely continuous) primitive $F$
  by an increasing sequence of Lipschitz functions $F_n$ uniformly converging to $F$, e.g. by setting
  \begin{displaymath}
    F_n(w)=\inf_{v\in [0,1]}F(v)+n|v-w|.
  \end{displaymath}
  Thus $F_n^{**}$ is an increasing sequence of Lipschitz functions satisfying $F_n^{**}(w)=F_n(w)$ at $w=0,1$, and
  pointwise converging to some lower semicontinuous convex function $G$ as $n\uparrow+\infty$ with
  \begin{equation}
    \label{eq:93}
    G(w)\le F^{**}(w)\le F(w)\quad\forall\, w\in [0,1].
  \end{equation}
  On the other hand, for $w=0,1$ we have $G(w)=\lim_{n\uparrow+\infty}F_n(w)=F(w)$ so that $F^{**}(w)=F(w)$.
  \eqref{eq:93} also yields
  \begin{displaymath}
  G(w)\le \liminf_{v\to w}F^{**}(v)\le \limsup_{v\to w}F^{**}(v)\le F(w)\quad\forall\, w\in [0,1]  
  \end{displaymath}
  so that $F^{**}$ is also continuous at $w=0,1$, where $G=F$.
\end{proof}
Let us now consider the set 
\begin{equation}
  \label{eq:32}
  \Lambda= \Bigl\{w\in [0,1]: (F - F^{**})(w) > 0 \Bigr\};
\end{equation}
since $\Lambda$ is open and does not contain $0$ and $1$, it is the disjoint union of a (at most countable) collection $\cO$ of
open intervals.
\begin{lemma} \label{F** linear}
  If $(a,b)\in \cO$ is a connected component of $\Lambda$ then for every $w=(1-\theta)a+\theta b$,
  $\theta\in [0,1]$
\begin{equation} \label{2-convex is const}
	F^{**}((1-\theta)a + \theta b) = (1-\theta )F(a) + \theta F(b)\quad
        \forall\, \theta\in [0,1],\qquad
        F(a)=F^{**}(a),\ F(b)=F^{**}(b).
\end{equation}
\end{lemma}
\begin{proof}
  Since $a,b\not\in \Lambda$ one has $F(a)=F^{**}(a)$ and $F(b)=F^{**}(b)$.
  Let $\bar w\in [a,b]$ a minimizer of the continuous function 
  $$
  w\mapsto F( w) - L(w),\quad L(w):= F(a)+(w-a)\frac{F(b)-F(a)}{b-a},
  $$
  so that $F(w)\ge F(\bar w)+L(w-\bar w)$ for every $w\in [a,b]$.
  The continuous function
  \begin{equation}
    \label{eq:33}
    G(w):=
    \begin{cases}
      F^{**}(w)&\text{if }w\not\in [a,b],\\
      \max\big(F^{**}(w),F(\bar w)+L(w-\bar w)\big)&\text{if }w\in [a,b],
    \end{cases}
  \end{equation}
  provides a convex lower bound of $F$ and therefore $G(w)\le F^{**}(w)$ for every $w\in [0,1]$.
  Since $G(\bar w)=F(\bar w)$ we deduce that $\bar w\not\in \Lambda$ and therefore $\bar w$ coincide
  with $a$ or $b$ and the inequality $G(w)\le F^{**}(w)$ yields
  $F^{**}((1-\theta)a + \theta b) \ge (1-\theta )F(a) + \theta F(b)$;
  the opposite inequality is a consequence of the convexity of $F^{**}$.
\end{proof} 
The next lemma contains the crucial inequality we need to characterize $\Proj\K$.
\begin{lemma}
  \label{le:monotonicity}
  Let $\psi\in C^1(\R)$ be a convex function. For every $f\in L^2(0,1)$ and $z\in \K$ with $g:=(F^{**})'$, if
  $(f-g)\psi'(z-g)\in L^1(0,1)$ we have
  \begin{equation}
    \label{eq:94}
    \int_0^1 \big(f(w)-g(w)\big)\psi'\big(z(w)-g(w)\big)\,\d w\le 0
    \le \int_0^1 \big(f(w)-g(w)\big)\psi'\big(g(w)-z(w)\big)\,\d w.
  \end{equation}
\end{lemma}
\begin{proof}
We decompose
$[0,1]$ in the disjoint union of the open intervals $(a,b)\in \cO$ covering $\Lambda$ (see \eqref{eq:32})
and of $[0,1]\setminus \Lambda$, where $F(w)=F^{**}(w)$, and therefore $f(w)=g(w)$
up to a $\Leb 1$-negligible set (recall that $F^{**}$ is locally Lipschitz).
In each $(a,b)\in \cO$ $F^{**}$ is linear, $g$ is constant, and
the function $w\mapsto \psi'(z(w)-g)$ is bounded and nondecreasing, thus its distributional derivative
is a nonnegative finite measure $\gamma_{a,b}$. Since $F=F^{**}$ in $\{a,b\}$, we have
\begin{align*}
  \int_0^1 &(f-g)\,\psi'(z-g) \,\d w
  = \int_\Lambda (f-g) \psi'(z-g) \,\d w + 
  \int_{[0,1]\setminus\Lambda} (f-g) \psi'(z-g) \,\d w 
  \\&= \sum_{(a,b)\in \cO}\int_a^b (f-g) \psi'(z-g)\, \d w 
  = -\sum_{(a,b)\in \cO} \int_a^b (F(w)-F^{**}(w))\,d\gamma_{a,b}(w)\le 0.
\end{align*}
The second inequality of \eqref{eq:94} can be simply obtained by
considering the convex function $\tilde\psi(r):=\psi(-r)$.
\end{proof}
\begin{proof}[End of the proof of Theorem 3.1]
  Concerning the projection in $L^2(0,1)$, by a standard approximation argument, it is not restrictive to assume $f \in L^\infty(0,1)$
so that $g\in L^\infty(0,1)$ too. Choosing $\psi(r):=\frac 12 r^2$ \eqref{eq:94} yields \eqref{eq:34}.

In order to prove \eqref{eq:96}, a standard approximation of $\psi$ by the increasing sequence
of its Moreau-Yosida approximations $\psi_n(r):=\min_{s\in \R} \psi(s)+\frac n2|s-r|^2$ shows that 
it is not restrictive to assume $\psi$ convex, $C^1,$ and at most quadratically growing as $|r|\to\infty$.
We can then apply the standard convexity inequality $\psi(s)-\psi(r)\ge \psi'(r)(s-r)$ and  Lemma \ref{le:monotonicity}
obtaining
\begin{align*}
  \int_\R &\Big(\psi\big(f-h\big)-\psi\big(\Proj\K(f)-\Proj\K(h)\big)\Big)\,\d w
  \\&\ge
  \int_\R \psi'\big(\Proj\K(f)-\Proj\K(h)\big)\Big(\big(f-\Proj\K(f)\big)-\big(h-\Proj\K(h)\big)\Big)\,\d w\topref{eq:94}\ge0.
\end{align*}
The first inequality of \eqref{eq:96} is a particular case of the second one, with $h=\Proj \K(h)=0$.
\end{proof}

The following result is a simple consequence of Theorem \ref{proj on K}.
Let us first introduce for a given $f\in L^2(0,1)$
the open set $\Omega_f\subset (0,1)$ where $f$ is locally constant
\begin{equation}
  \label{eq:2}
  \Omega_f:=\big\{w\in (0,1):f\text{ is essentially constant in a neighborhood of }w\big\}.
\end{equation}
Equivalently $\Omega_f$ is the complement of the support of the distributional derivative of $f$.
\begin{corollary}
  \label{cor:omega}
  Let $f\in L^2(0,1)$ and $g=\Proj \K(f)$. Then
  \begin{equation}
    \label{eq:69}
    \Omega_f\subset \Omega_g.
  \end{equation}
\end{corollary}
\begin{proof}
  Notice that $\Lambda\subset \Omega_g$ ($\Lambda$ has been defined by
  \eqref{eq:32}); if $w\in \Omega_f\setminus \Lambda$ then $F(w)=F^{**}(w)$, so that
  any linear part of the graph of $F$ in an open interval containing $w$ should locally coincide with $F^{**}$; it follows that
  $F^{**}=F$ in a neighborhood of $w$ so that $w\in \Omega_g$.
\end{proof}
\begin{definition}[The polar cone and the subdifferential of the indicator function $I_\K$]
  \label{def:polar}
  We denote by $\K^\circ$ the polar cone of $\K$, defined by
  \begin{equation}
    \label{eq:42}
    f\in \K^\circ\quad \Longleftrightarrow\quad
    \Scalar fz\le 0\quad \forall\, z\in \K\quad
    \Longleftrightarrow\quad
    \Proj \K(f)=0.
  \end{equation}
  The subdifferential $\partial I_\K(g)$ of the indicator function of $\K$ (see \eqref{eq:22})
  at some
  function $g\in \K$ is the subset of $L^2(0,1)$ characterized by
  \begin{equation}
    \label{eq:117}
    \xi \in \partial I_\K(g) \quad\Longleftrightarrow\quad
    \Scalar\xi{z-g}\le 0\quad\forall\, z\in \K.
  \end{equation}
\end{definition}
\begin{remark}
    \label{rem:polar}
    $\K^\circ$ and $\partial I_\K$ are clearly linked by $\K^\circ=\partial I_\K(0)$ and
  \begin{equation}
    \label{eq:43}
    \xi \in \partial I_\K(g)\quad \Longleftrightarrow\quad
    \xi\in \K^\circ,\quad \Scalar \xi g=0.
  \end{equation}
  $\K^\circ$ provides an equivalent reformulation of \eqref{eq:34}, since
  \begin{equation}
    \label{eq:50}
    g=\Proj \K(f)\quad\Longleftrightarrow\quad g\in K,\quad
    f-g\in \K^\circ,\quad \Scalar{f-g}g=0
    \quad
    \Longleftrightarrow\quad
    f-g\in \partial I_\K(g).
  \end{equation}
\end{remark}
If $\Omega$ is an open subset of $(0,1)$, we denote by $\PC \Omega$ the convex cone
\begin{equation}
  \label{eq:56}
  \PC\Omega:=\Big\{F\in C^0([0,1]):F\ge0\text{ in } [0,1],\quad
  F=0\text{ in }[0,1]\setminus\Omega\Big\}.
\end{equation}
We can give a useful characterization of $\K^\circ$ in term of the cone $\PositiveCone:=\PC{(0,1)}$.
\begin{proposition}[A characterization of the polar cone $\K^\circ$]
  \label{prop:polar}
  A function $f$ belongs to the polar cone $\K^\circ$ if and only if its primitive $F(w):=\int_0^w f(s)\,\d s$
  belongs to $\PositiveCone$.
\end{proposition}
\begin{proof}
  If $F\in \PositiveCone$ then one easily gets for every $z\in \K\cap C^1([0,1])$
  \begin{equation}
    \label{eq:57}
    \Scalar fz=\int_0^1 F'(w)\,z(w)\,\d w=-\int_0^1 F(w)\,z'(w)\,\d w\le 0,
  \end{equation}
  since $F,z'\ge 0$, $F(0)=F(1)=0$.

  Let us now assume that $f\in \K^\circ$; for every continuous and nonnegative
  function $z\ge0$ and $c\in \R$ with $Z(w)=\int_0^w z(s)\, \d s-c$, since $Z\in \K$ we have
  \begin{displaymath}
    0\ge \Scalar fZ=\int_0^1 f(w)\, Z(w)\,\d w=-\int_0^1 F(w) \,z(w)\, \d w+F(1)(Z(1)-c)
  \end{displaymath}
  Since $c,z$ are arbitrary, we conclude that $F\in \PositiveCone$.
\end{proof}
The last result of this section concerns a precise characterization of
$\partial I_\K$. Let us first define for $f\in L^2(0,1)$ the closed subspace $\cH_f\subset L^2(0,1)$
defined as
\begin{equation}
  \label{eq:3}
  \cH_f:=\big\{h\in L^2(0,1): h\text{ is essentially constant in each connected component of }\Omega_f\big\}.
\end{equation}
We denote by $\ProjH f$ the orthogonal $L^2$-projection on $\cH_f$. It is easy to check that
\begin{equation}
  \label{eq:58}
  \begin{aligned}
    \ProjH g(f)&=f\text{ a.e.\ in }(0,1)\setminus \Omega_g,\\
    \ProjH g(f)&\equiv \media_\alpha^\beta f(w)\,\d w\text{ a.e.\ in every
      connected component }(\alpha,\beta)\subset \Omega_g.
  \end{aligned}
\end{equation}
Moreover, denoting by $F$ the primitive function of $f$,
\begin{equation}
  \label{eq:60}
  \text{if }F\in \PC{\Omega_g}\text{ then $f$ is orthogonal to }\cH_g,
\end{equation}
since $f=F'$ vanishes a.e.\ outside $\Omega_g$ and for every connected component $(\alpha,\beta)$ of
$\Omega_g$ we have $\int_\alpha^\beta f(w)\, \d w=F(\beta)-F(\alpha)=0$.
\begin{theorem}[The subdifferential of $I_\K$]
  \label{thm:minimal}
  Let $g\in \K$, $\xi\in L^2(0,1)$, and $\Xi(w):=\int_0^w \xi(s)\,\d s$.
  Then we have
  \begin{equation}
    \label{eq:52}
    \xi\in \partial I_\K(g)
    \quad\Longleftrightarrow\quad
    \Xi\in \PC{\Omega_g}.
  \end{equation}
  In particular,
  \begin{equation}
    \label{eq:59}
    \text{ if $\xi\in \partial I_\K(g)$ then}
    \quad
    \left\{
    \begin{aligned}
      &\xi=0\text{ a.e.\ in } [0,1]\setminus\Omega_g,\\
      &\int_\alpha^\beta \xi(w)\,\d w=0\quad\text{for every connected component $(\alpha,\beta)$ of $\Omega_g$,}
      \end{aligned}
      \right.
  \end{equation}
  so that $\xi$ is orthogonal to $\cH_g$ and we have by \eqref{eq:50} and \eqref{eq:69}
  \begin{equation}
    \label{eq:82}
    g=\Proj\K(f)\quad\Rightarrow\quad g=\Proj{\cH_g}(f),\quad
    \cH_g\subset \cH_f.
  \end{equation}
\end{theorem}
\begin{proof}
  The left implication in \eqref{eq:52} is immediate, since $\Xi\in \PC{\Omega_g}$ implies
  $\Xi\in \PositiveCone$ and therefore $\xi\in \K^\circ$ by Proposition \ref{prop:polar};
  moreover, $\xi$ is orthogonal to $\cH_g$ by \eqref{eq:60} and therefore it is also orthogonal
  to $g\in \cH_g$, so that $\xi \in \partial I_\K(g)$ by \eqref{eq:43}.

  Conversely, if $\xi \in \partial I_\K(g)$, then $\Xi\in \PositiveCone$ by \eqref{eq:43}
  and Proposition \ref{prop:polar}. Moreover, denoting by $\gamma=g'$ the nonnegative Radon measure
  associated to the distributional derivative of $g$ in $(0,1)$,
  the next Lemma \ref{le:Hardy} yields
  \begin{equation}
    \label{eq:61}
    0\topref{eq:43}=\int_0^1 \xi(w)g(w)\,\d w\topref{eq:62}=-\int_0^1 \Xi(w)\,\d\gamma(w),
  \end{equation}
  which shows that $\Xi(w)=0$ on the support of $\gamma$ and yields $\Xi\in \PC{\Omega_g}$.
\end{proof}
\begin{lemma}
  \label{le:Hardy}
  Let $g\in \K$ and $\xi\in \K^\circ$ with (nonnegative) primitive $\Xi\in \PositiveCone$.
  If $\gamma=g'$ is the nonnegative Radon measure associate to the distributional derivative of $g$
  in $(0,1)$ then $\Xi\in L^1(\gamma)$ and
  \begin{equation}
    \label{eq:62}
    \int_0^1 g(w)\xi(w)\,\d w=-\int_0^1 \Xi(w)\,\d\gamma(w).
  \end{equation}
\end{lemma}
\begin{proof}
  Since $\gamma$ is a nonnegatative Radon measure in $(0,1)$ but not necessarily finite, 
  we need an approximation argument to justify \eqref{eq:62}. Let $\varphi_n\in C^\infty_0(0,1)$ be
  an increasing sequence of nonnegative functions such that $\lim_{n\up+\infty}\varphi_n(w)=1$, $|\varphi_n'|\le 2 n$,
  and $\varphi_n(w)\equiv 1$ for $1/n\le w\le 1-1/n$.
  We have
  \begin{subequations}
    \begin{align}
      \label{eq:64}
      \int_0^1 g\xi\varphi_n\,\d w&=-\int_0^1 \Xi\varphi_n\,\d\gamma-
      \int_0^1 \Xi g\varphi_n'\,\d w \\&=\label{eq:63} -\int_0^1
      \Xi\varphi_n\,\d\gamma- \int_0^{1/n} \Xi g\varphi_n'\,\d w-
      \int_{1-1/n}^1 \Xi g\varphi_n'\,\d w.
    \end{align}
  \end{subequations}
  Applying Hardy inequality, we get
  \begin{displaymath}
    \left|\int_0^{1/n} \Xi g\varphi_n'\,\d w\right|\le 2n \|w^{-1}\Xi\|_{L^2(0,1/n)}\|wg\|_{L^2(0,1/n)}\le
    2C \|\xi\|_{L^2(0,1)}\|g\|_{L^2(0,1/n)}
  \end{displaymath}
  so that the integral vanishes as $n\uparrow+\infty$. A similar argument holds for the last integral of \eqref{eq:63}.
  Passing to the limit in (\ref{eq:64},b) as $n\uparrow+\infty$ and using Lebesgue dominated (being
  $g\xi\in L^1(0,1)$) or monotone (being $\Xi\ge0$ and $\varphi_n$ increasing) convergence theorem, we conclude.
\end{proof}
The last Lemma of this section provides a useful example concerning a class of elements in $\partial I_\K(g)$.
\begin{lemma}[An example of minimal selection in $\partial I_\K$]
  \label{le:antimonotone}
  If $g,h\in \K$, then
  \begin{equation}
    \label{eq:65}
    \xi_h:=\ProjH g(h)-h\in \partial I_\K(g).
  \end{equation}
  Moreover, 
  \begin{equation}
    \label{eq:66}
    \|z-h-\xi_h\|_{L^2(0,1)}\le \|z-h-\xi\|_{L^2(0,1)}\quad\forall\, \xi\in \partial I_\K(g), \quad z\in \cH_g.
  \end{equation}
  In particular,
  \begin{equation}
    \label{eq:79}
    \text{if $z\in \cH_g$ then}\quad
    \|z\|_{L^2(0,1)}\le \|z-\xi\|_{L^2(0,1)}\quad\forall\, \xi\in \partial I_{K}(g).
  \end{equation}
\end{lemma}
\begin{proof}
  Since $h-\ProjH g(h)$ is orthogonal to $\cH_g$ (thus in particular to $g$), by \eqref{eq:43} we have to check that
  $\xi_h\in \K^\circ$, by applying Proposition \ref{prop:polar}. By \eqref{eq:58}, $\xi_h=0$ a.e.\ in $(0,1)\setminus \Omega_g$, so that
  the primitive $\Xi_h$ of $\xi_h$ satisfies
  \begin{displaymath}
    \Xi_h(w)=\int_{\Omega_g\cap (0,w)}\xi_h(s)\,\d s.
  \end{displaymath}
  The thesis then follow if we show that for every connected component
  $(\alpha,\beta)$ of $\Omega_g$ we have $\Xi_h(\alpha)=\Xi_h(\beta)=0=\min_{[\alpha,\beta]} \Xi_h$.
  Since the characteristic function $\nchi_{(0,\alpha)}$ of $(0,\alpha)$ belongs to $\cH_g$, we have
  \begin{displaymath}
    \Xi_h(\alpha)=\int_0^\alpha \xi_h(w)\,\d w=\Scalar{h-\ProjH g(h)}{\nchi_{(0,\alpha)}}=0.
  \end{displaymath}
  A similar argument shows that $\Xi_h(\beta)=0$. Moreover, for $w\in (\alpha,\beta)$ we have
  \begin{displaymath}
    \Xi_h(w)=\int_\alpha^w \xi_h(s)\,\d s\topref{eq:58}=(w-\alpha)\media_\alpha^\beta h(w)\,\d w-\int_\alpha^w h(w)\,\d w,
  \end{displaymath}
  which shows that $\Xi_h$ is concave, and therefore nonnegative in $(\alpha,\beta)$.

  \eqref{eq:66} follows immediately by observing that
  $\xi,\xi_h\in(\cH_g)^\perp$ and
  $z-h-\xi_h$ belongs to $\cH_g$ and therefore it is the orthogonal projection of
  $z-h$ onto $(\cH_g)^\perp$.
\end{proof}

\section{The Lagrangian formulation of the discrete sticky particle system}
\label{sec:dSPS}
In this section we shall show that
the discrete sticky particle system satisfies the three characterizations
of Theorem \ref{thm:L2flow} and we prove Theorem \ref{thm:main1}.
\begin{notation}
  \label{not:sticky}
  Let us recapitulate our basic notation and definitions
  \begin{enumerate}
  \item $P_i(t)=(m_i,x_i(t),v_i(t))$, $i\in I=\{1,\cdots,N\}$, $t\ge
    0$, is a solution of the discrete sticky particle system;
  \item the positions of the particles are ordered: $x_1(t)\le
    x_2(t)\le \ldots\le x_N(t);$
  \item the sets $J_i(t):=\{j\in I: x_j(t)=x_i(t)\}$ are nondecreasing
    with respect to time.  They correspond to a single particle of
    mass $\sum_{j\in J_i(t)}m_j$.
  \item At each time $t$ we pick up the collection of minimal indexes
    \begin{displaymath}
      I(t):=\big\{\min J_i(t):i=1,\dots,N\big\}=\big\{i_1(t)<\ldots < i_{N(t)}\big\}\subset I,
    \end{displaymath}
    so that each $J_i(t)$ is of the form $\{j\in I: i_k(t)\le j<i_{k+1}(t)\}$ for some $k$ and
    $\big(J_i(t)\big)_{i\in I(t)}$ is a partition of $I$.
  \item We denote by $0<t_1<t_2<\ldots t_h<\ldots<t_{H-1}$ the
    (finite) sequence of times at which the cardinality of some
    $J_i(t)$ has an increasing jump; setting $t_0=0$ and
    $t_H=+\infty$, $\{[t_h,t_{h+1})\}_{h=0}^{H}$ is the associate
    partition of the positive real line with step sizes $\delta_h:=t_{h}-t_{h-1}.$
  \item The functions $x_i$ are continuous and piecewise linear on
    each interval $[t_h,t_{h+1})$, with piecewise constant, right continuous derivatives $v_i(t)$
    satisfying \eqref{eq:1}.
    Each set $J_i(t)$ and $I(t)$ is also constant
    in each interval $[t_h,t_{h+1})$.  
  \end{enumerate}
\end{notation}
Let $\rho_t=\sum_{i\in I}m_i\delta_{x_i(t)}$ be the measure induced by the discrete
sticky-particle system. In order to write explicitly the function $\X(t):=\X_{\rho_t}$ we consider the subdivision of $[0,1]$
given by
\begin{equation}
  \label{eq:35}
  w_0=0<w_1<\ldots<w_N=1,\quad w_{i}=w_{i-1}+m_i=\sum_{j=1}^i m_j,\qquad i\in I.
\end{equation}
We also set
\begin{equation}
  \label{eq:36}
  W_i:=[w_{i-1},w_i),\quad
  W_i(t)=\bigcup_{j\in J_i(t)} W_j,\qquad i\in I,
\end{equation}
and we notice that
\begin{equation}
  \label{eq:37}
  \X(t)=\sum_{i=1}^N x_i(t)\Char_{W_i},\quad
  \Rightdert\X(t)=V(t)=\sum_{i=1}^N v_i(t)\Char_{W_i}.
\end{equation}
The main result of this section is
\begin{theorem}[Lagrangian formulation of the discrete SPS]
  \label{thm:discrete_Lagrangian}
  The couple $(X,V)$ defined by \eqref{eq:37} satisfies the equations
  \emph{(L.I, II, III)} and the properties \emph{(L.a,b,c)} of Theorem \ref{thm:L2flow}.
  In particular, it defines a semigroup $\sfS_t$ in the discrete subspace
  \begin{equation}
    \label{eq:115}
    \hat \cX:=\Big\{(X,V)\in \Lspace p:
    \X=\sum_{i=1}^N x_i\Char_{W_i} \
    \text{for a finite interval partition $(W_i)_{i=1}^N$ of $[0,1)$}\Big\}
  \end{equation}
\end{theorem}
We split the \emph{proof} in various steps.

The collection $\big(W_i(t)\big)_{i\in I(t)}$ is a partition of $[0,1).$
In $L^2(0,1)$ we introduce the decreasing family of finite dimensional spaces $\cH(t)$ whose elements are
piecewise constant on each interval $W_i(t)$, $i\in I(t)$.
Notice that, by the very definitions of $\Omega_{X(t)}$ and $\cH_{X(t)}$ \eqref{eq:2} and
\eqref{eq:3}
\begin{equation}
  \label{eq:75}
  \Omega_{X(t)}=(0,1)\setminus \{w_i:i\in I(t)\},\quad
  \cH(t)=\cH_{X(t)}.
\end{equation}
Besides \eqref{eq:37}, the crucial features describing the evolution of $X(t)$ are
\begin{equation}
  \label{eq:45}
  \X(t)\in \K\cap \cH(t),\quad V(t)\in \cH(t),\quad
  \cH(t)=\cH_h,\ V(t)=V_h\text{ if }t\in [t_h,t_{h+1}),
\end{equation}
and the update rule for the velocity \eqref{eq:1}: $V(t_{h})$ is constant
in each interval $W_i(t_h)=\cup_{j\in J_i(t_h)}W_j$ and its value is given by
\begin{align*}
    V(t_h+)\restr{W_i(t_h)}&=\frac{\sum_{j\in J_i(t_h)} m_j v_j(t_{h-1})}{\sum_{j\in J_i(t_{h})} m_i}=
    \big(\Leb 1(W_i(t_h)\big)^{-1}\int_{W_i(t_h)} V(t_{h-1})\,\d w, 
\end{align*}
so that by \eqref{eq:58}
\begin{equation}
  \label{eq:76}
  V_{h}=\Proj{\cH_{h}}(V_{h-1})=\Proj{\cH_{h}}(V_0)\quad\text{since}\quad
  \cH_0\supset \cH_1\supset \cH_2\supset\cdots\cH_h,
\end{equation}
which yields \eqref{eq:74} and \eqref{eq:108}.
The next lemma shows \eqref{eq:26}.
\begin{lemma}
  \label{le:representation}
  Let $\tilde \X(t):=\X_0+tV_0$ be associated to the free system
  $\tilde P_i=(m_i,\x_i,\v_i)$ given by $\x_i(t)=x_i(0)+tv_i(0),\ \v_i(t)\equiv \v_i= v_i(0)$.
  Then
  \begin{equation}
    \label{eq:46}
    \X(t)=\Proj{\cH(t)}(\tilde X(t))=\Proj{\cH(t)}(X_0+tV_0),
  \end{equation}
  \begin{equation}
  \label{eq:40}
  t\Rightdert \X(t)=tV(t)=\X(t)-X_0-\Xi(t)\quad\text{for}
  \quad\Xi(t):= -\X_0+\Proj{\cH(t)}(\X_0)
  \in \partial I_\K(\X(t)).
  \end{equation}
\end{lemma}
\begin{proof}
  Suppose that $t\in [t_h,t_{h+1})$; since $X(t)\in \cH(t)=\cH_h\subset \cH(r)$ and
  $V(r)=\Proj{\cH(r)}(V_0)$
  for $0\le r\le t$ by \eqref{eq:76}, we have by the linearity of $\Proj{\cH(r)}$
  \begin{align*}
    X(t)&=X_0+\int_0^t V(r)\,\d r\topref{eq:45}=
    \Proj{\cH(t)}\Big(X_0+\int_0^t V(r)\,\d r\Big)\topref{eq:76}=
    \Proj{\cH(t)}(X_0)+\int_0^t \Proj{\cH(t)}(\Proj{\cH(r)}(V_0))\,\d r
    \\&=
    \Proj{\cH(t)}(X_0)+\int_0^t \Proj{\cH(t)}(V_0)\,\d r=
    \Proj{\cH(t)}(X_0+\int_0^t V_0\,\d r)=\Proj{\cH(t)}(\tilde X(t)).
  \end{align*}
  From of \eqref{eq:46} we have
  \begin{align*}
    t\Rightdert \X(t)=tV(t)\topref{eq:76}=\Proj{\cH(t)}{(tV_0)}
    \topref{eq:46}=X(t)-
    \Proj{\cH(t)}{(X_0)}=
    X(t)-X_0-\Xi(t),
  \end{align*}
  where $\Xi(t)=\Proj{\cH(t)}(X_0)-X_0$; since $X_0\in \K$ and
  $\cH(t)=\cH_{X(t)}$,
  by Lemma \ref{le:antimonotone} we conclude that $\Xi(t)\in \partial I_\K(X(t))$.
\end{proof}
We conclude now the proof of \eqref{eq:44} and \eqref{eq:49}; notice that \eqref{eq:49bis}
follows directly by \eqref{eq:49} and \eqref{eq:74} via the semigroup property
of $\SS_t$ in $\DXspace(\R)$.
\newcommand{\Ui}{U}.
\begin{lemma}
  \label{le:2}
  Under the same notation and assumptions as before, we have
  \begin{equation}
    \label{eq:47}
    \Ui(t):=V_0-V(t)=V_0-\Proj{\cH(t)}(V_0)\in \partial I_\K(X(t)),\quad
    X(t)=\Proj\K(X_0+tV_0).
  \end{equation}
\end{lemma}
\begin{proof}
  Since $\Scalar{\Ui(t)}{X(t)}=0$ (being $X(t)\in \cH(t)$), the first inclusion of \eqref{eq:47} is equivalent to
\begin{equation}
  \label{eq:51}
  \Ui(t)= V_0-V(t)\in \K^\circ\quad\forall\, t\ge0,
\end{equation}
by \eqref{eq:43}.
It is not restrictive to assume that $t=t_h$ and $V(t)=V_h$ for some $h\in\{1,\ldots,H-1\}.$
Since $\K^\circ$ is a cone and $V_0-V_h$ can be decomposed into the sum
\begin{equation}
  \label{eq:48}
  V_0-V_h=\sum_{k=0}^{h-1}(V_k-V_{k+1})
\end{equation}
it is sufficient to prove that $V_k-V_{k+1}\in \K^\circ$ or, equivalently, that $\delta_{k+1}(V_k-V_{k+1})\in \K^\circ$.
Since $V_{k+1}=\Proj{\cH_{k+1}}(V_k)$ we obtain
\begin{align*}
  \delta_{k+1}(V_k-V_{k+1})&=\delta_{k+1}V_k-\Proj{\cH_{k+1}}(\delta_{k+1}V_k)=
  (\X_{k+1}-\X_k)-\Proj{\cH_{k+1}}(\X_{k+1}-\X_k)\\&=
  \X_{k+1}-\Proj{\cH_{k+1}}(X_{k+1})+\Proj{\cH_{k+1}}(\X_k)-X_k=
  \Proj{\cH_{k+1}}(\X_k)-X_k\in \K^\circ
\end{align*}
by Lemma \ref{le:antimonotone}.

The second identity of \eqref{eq:47} follows now by a similar argument, by checking the conditions of \eqref{eq:50}.
Since $X(t)\in \K$ and $\Scalar{\tilde X(t)-\X(t)}{\X(t)}=0$ by
\eqref{eq:46}, it is sufficient to show that
$
  \tilde X(t)-\X(t)\in \K^\circ.
$ 
On the other hand
\begin{displaymath}
  \tilde \X(t)-\X(t)=tV_0-\int_0^t V(r)\,\d r=
  \int_0^t\big( V_0-V(r)\big)\,\d r=\int_0^t \Ui(r)\,\d r
\end{displaymath}
and \eqref{eq:51} shows that $\Ui(r)\in \K^\circ$ for every $r\ge0$.
Being $\K^\circ$ a cone, we conclude.
\end{proof}
\begin{proof}[\textbf{Proof of Theorem \ref{thm:main1}}]
Let us now consider two discrete Lagrangian solutions $(X^\ell(t),V^\ell(t))=\sfS_t(X^\ell_0,V^\ell_0)\in \hat\cX$, $\ell=1,2$.
\eqref{eq:96}, \eqref{eq:102}, and 
\eqref{eq:49} immediately yield the estimates
\begin{equation}
  \label{eq:104}
  \int_0^1 \psi\big(X^1(t)-X^2(t)\big)\,\d w\le
  \int_0^1 \psi\big(X^1_0-X^2_0+t(V^1_0-V^2_0)\big)\,\d w
\end{equation}
\begin{equation}
  \label{eq:77}
  \|X^1(t)-X^2(t)\|_{L^p(0,1)}\le \|X^1_0-X^2_0\|_{L^p(0,1)}+t\|V^1_0-V^2_0\|_{L^p(0,1)},
\end{equation}
which are equivalent to \eqref{eq:98} and \eqref{eq:16}.
\eqref{eq:44} yields (\cite[Theorem 3]{Savare93}, \cite[Theorem 1.2]{Savare96})
\begin{equation}
  \label{eq:78}
  \int_0^t \|V^1-V^2\|^2\,\d r\le C(1+t) \Big(\sum_{\ell=1,2}\|X^\ell_0\|
  +\|V^\ell_0\|
  \Big)
  \Big(\|X^1_0-X^2_0\|+\|V^1_0-V^2_0\|\Big),
\end{equation}
which is equivalent to \eqref{eq:17}.
\end{proof}
\section{Stability and uniqueness of Lagrangian solutions}
\label{sec:lagrangian}
Our first result concerns the stability of Lagrangian solutions
to (L.I, II, III) of
Theorem \ref{thm:L2flow} (in particular it applies to those obtained by the
discrete SPS in $\hat\cX$).
\begin{lemma}
  \label{le:Lagrangian_stability}
  Let $X^n,V^n:=\Rightdert X^n$ curves satisfying
  all the equations \emph{(L.I, II, III)} and the properties
  \emph{(L.a,b,c)} stated in Theorem \ref{thm:L2flow} with respect to
  initial data $X^n_0,V^n_0=v^n_0(X^n_0)$ converging to
  $X_0,V_0=v_0(X_0)$ in $L^p(0,1)$, $p\ge 2$.
  \begin{enumerate}[(a)]
  \item 
    $X^n(t)$ converges to $X(t)$ in $L^p(0,1)$, uniformly in each
    compact interval; $X$ is Lipschitz continuous with values in $L^p(0,1)$.
  \item The Lipschitz curve $X$ is 
    right-differentiable at each point $t$, with right-continuous derivative $V(t)$, and
    it satisfies \emph{(L.I, II, III)} and \emph{(L.a,b,c)} of Theorem
    \ref{thm:L2flow}.
    \item $V^n$ strongly converges to $V$ in
    $L^2(0,T;L^2(0,1))$ for every $T>0$.
  \item The
    curve $X$ is differentiable in $L^p(0,1)$ and $V$ is continuous at each point of
    $(0,+\infty)\setminus\cT$,  where $\cT$ is the jump
    set of the nonincreasing map $t\mapsto \|V(t)\|_{L^2(0,1)}$.
    \item If $\bar V$ is any weak accumulation point of $V^n(t)$ in $L^p(0,1)$, then
    $\Proj{\cH_{X(t)}}(\bar V)=V(t)$.      
  \item $V^n(t)\to V(t)$ in $L^p(0,1)$
    for every $t\in [0,+\infty)\setminus \cT$.
  \end{enumerate}
\end{lemma}
\begin{proof}
  \emph{(a)} is an immediate consequence of
  \eqref{eq:49} and \eqref{eq:102}, which also show that
  $X^n$ is uniformly Lipschitz continuous with values in $L^p(0,1)$ and Lipschitz constant
  bounded by $\|V^n_0\|_{L^p(0,1)}$. The convergence is therefore uniform in each compact interval
  and the limit function $X$ satisfies the same Lipschitz bound with constant $\|V_0\|_{L^p(0,1)}$.

  \emph{(b,c)} Standard stability results for gradient flows in Hilbert spaces \cite{Brezis73} show that $X$ solves
  \eqref{eq:44} and \eqref{eq:26}; in particular $X$ is right differentiable in $L^2(0,1)$ at each $t\ge0$,
  with $L^2(0,1)$-right derivative $V(t)$ which is right-continuous.
  \eqref{eq:78} shows that $V$ is the limit of $V_n$ in $L^2(0,T;L^2(0,1))$ for every
  $T>0$ (this proves point \emph{(c)}): in particular, up to the extraction of a suitable subsequence $n_k$, we can find an $\Leb 1$-negligible set
  $N\subset (0,+\infty)$ such that $V_{n_k}(t)\to V(t)$ in $L^2(0,1)$ for every $t\in [0,+\infty)\setminus N$ as $k\uparrow+\infty$.
  Passing to the limit in \eqref{eq:49bis} and in \eqref{eq:44} we obtain that
  \begin{equation}
    \label{eq:83}
    X(t)=\Proj\K(X(s)+(t-s)V(s)),\qquad
    \Rightdert X(t)=V(t)\in -\partial I_\K(X(t))+V(s)
  \end{equation}
  for every $s\in [0,+\infty)\setminus N$ and $t\ge s$. Since $V$ is right continuous, \eqref{eq:83}
  eventually holds for every $0\le s\le t$.

  The projection formula of \eqref{eq:83} shows that
  \begin{equation}
    \label{eq:68}
    \|X(t+h)-X(t)\|_{L^p(0,1)}\le h\|V(t)\|_{L^p(0,1)}\le h\|V(s)\|_{L^p(0,1)}\quad \forall\, 0\le s\le t,\quad h\ge 0,
  \end{equation}
  and, more generally,
  \begin{equation}
    \label{eq:106}
    \int_0^1\psi\big(h^{-1}(X(t+h)-X(t))\big)\,\d w\le
    \int_0^1 \psi(V_s)\,\d w\le \int_0^1 \psi(V_0)\,\d w\quad \forall\, 0\le s\le t,\ h\ge0
  \end{equation}
  for every convex nonnegative function $\psi:\R\to \R$.
  \eqref{eq:68} and the right-differentiability of $X$ in $L^2(0,1)$
  yields that $V(t)$ is also the right derivative of $X$ in $L^p(0,1)$, its $L^p$ norm is not increasing,
  and by \eqref{eq:106} the family $V_s$ is uniformly $p$-integrable (by Dunford-Pettis criterion,
  it is sufficient to choose
  a convex function $\psi$ with $\psi(r)/|r|^p\to+\infty$ as $|r|\to+\infty$ and $\psi\circ V_0\in L^1(0,1)$,
  see e.g.\ \cite[Lemma 3.7]{Rossi-Savare03})

  From \eqref{eq:26} we deduce that $ t V(t)=X(t)-X_0-\Xi(t)$ where $\Xi(t)$ is characterized by
\begin{equation}
  \label{eq:84}
  \Xi(t)\in\partial I_\K(X(t)),\qquad \|X(t)-X_0-\Xi(t)\|\le
  \|X(t)-X_0-\xi\|\quad \forall\, \xi\in \partial I_\K(X(t)).
\end{equation}
Applying Lemma \ref{le:antimonotone} with $g:=X(t)$ and $h:=X_0$, we obtain
$\Xi(t)=\Proj{\cH_{X(t)}}(X_0)-X_0$ and therefore
\begin{equation}
  \label{eq:85}
  t V(t)=X(t)-\Proj{\cH_{X(t)}}(X_0),\quad
  V(t)\in \cH_{X(t)}.
\end{equation}
It follows by \eqref{eq:82} that
$\cH_{X(s)}\supset \cH_{X(t)}$ if $0\le s\le t$; moreover, by \eqref{eq:116},
there exists a Borel map $v_t\in L^p_{\rho_t}(\R)$ such that
\begin{equation}
  \label{eq:86}
  V(t)=v_t\circ X(t),\quad
  V(t)=\Proj{\cH_{X(t)}}(V(s))\quad\forall\,0\le s\le t,
\end{equation}
where the last identity follows by the fact that $V(t)$ belongs to $\cH_{X(t)}$ and
$\partial I_\K(X(t))$ is orthogonal to $\cH_{X(t)}$.

  \emph{(d)} Let $\cT$ be the jump set of the $L^2$-norm of $V(t)$; we show that $V$ is left-continuous
  at every $\bar t\in (0,+\infty)\setminus \cT$ (this also yields
  the left-differentiability of $X$ at $\bar t$).
  \eqref{eq:44} provides the minimal selection characterization of $V$
  \begin{equation}
    \label{eq:95}
    V(t)\in V_0-\partial I_\K(X(t)), \quad
    \|V(t)\|_{L^2(0,1)}\le \|V_0+ \xi\|_{L^2(0,1)}\quad\forall\, \xi\in \partial I_\K(X(t))\quad\forall\, t\ge0.
  \end{equation}
  Take an arbitrary increasing sequence $t_n\uparrow \bar t$ such that $V(t_n)\weakto \bar V$ 
  in $L^p(0,1)$. Since the graph of $\partial I_\K$ is strongly-weakly closed in $L^2(0,1)$,
  we have $\bar V\in V_0-\partial I_K(X(\bar t))$.
  Passing to the limit in \eqref{eq:95} we obtain
  \begin{equation}
    \label{eq:103}
    \|\bar V\|_{L^2(0,1)}\le \lim_{n\to\infty}\|V(t_n)\|_{L^2(0,1)}=
    \|V(\bar t)\|_{L^2(0,1)}\le  \|V_0+\bar \xi\|_{L^2(0,1)}\quad\forall\, \bar \xi\in \partial I_\K(X(\bar t)).
  \end{equation}
  Since $\partial I_K(X(\bar t))$ is a closed convex set, it follows that $\bar V=V(\bar t)$
  and the convergence is strong in $L^2(0,1)$ and therefore also in $L^p(0,1)$, since
  $V(t_n)$ is uniformly $p$-integrable.
    
  \emph{(e)} Let $n_k$ be an arbitrary subsequence such that
  $V^{n_k}(t)\weakto \bar V$ in $L^p(0,1)$.
  Passing to the limit in the inclusion $V^n(t)\in V^n_0-\partial I_\K(X^n(t))$
  we obtain $\bar V\in V_0-\partial I_\K(X(t))$. By Theorem \ref{thm:minimal}
  any element in $\partial I_\K(X(t))$ is orthogonal to $\cH_{X(t)}$ so that 
  $\Proj{\cH_{X(t)}}(\bar V)=\Proj{\cH_{X(t)}}(V_0)\topref{eq:74}=V(t)$.
    
    \emph{(f)} Let now $t\in (0,+\infty)\setminus \cT$ and let $n_k,\bar V$ be as
    in the previous point \emph{(e)}.
    Up to the extraction of a further
    subsequence (still denoted by $n_k$), there exists
  a dense set $S\subset (0,+\infty)$ such that $V^{n_k}(s)\to V(s)$ for every $s\in S$, so that
  \begin{displaymath}
    \|\bar V\|_{L^2(0,1)}\le \limsup_{k\uparrow+\infty}\|V^{n_k}(t)\|_{L^2(0,1)}\le
    \limsup_{k\uparrow +\infty}\|V^{n_k}(s)\|_{L^2(0,1)}=\|V(s)\|_{L^2(0,1)}\quad
    \forall\, s\in S,\ s<t.
  \end{displaymath}
  Since $t$ is a continuity point for $V$ we obtain by \eqref{eq:95}
   \begin{equation}
    \label{eq:95bis}
     \|\bar V\|_{L^2(0,1)}\le \|V(t)\|_{L^2(0,1)}
     \le \|V_0- \xi\|_{L^2(0,1)}\quad\forall\, \xi\in \partial I_\K(X(t)),
  \end{equation}
  which yields $\bar V=V(t)$, $\limsup_{k\uparrow+\infty}\|V^{n_k}(t)\|_{L^2(0,1)}\le
  \|V(t)\|_{L^2(0,1)}$, and the strong convergence of $V^n(t)$ to $V(t)$ in
  $L^2(0,1)$. The strong convergence in $L^p(0,1)$ follows by
  the uniform $p$-integrability estimate \eqref{eq:106}.
\end{proof}
\begin{corollary}[Existence of the Lagrangian semigroup]
  \label{cor:existence}
  For every initial data $(X_0,V_0)\in \Lspace 2$
  there exists a unique Lipschitz curve $X$ in $L^2(0,1)$ satisfying
  the equations \emph{(L.I, II, III)} and the properties \emph{(L.a,b,c)} stated in
  Theorem \ref{thm:L2flow}. Setting $V(t):=\Rightdert X(t)$, the map
  $\sfS_t:(X_0,V_0)\mapsto (X(t),V(t))$ defines a right-continuous semigroup in each
  space $\Lspace p$, $p\ge 2$.  
\end{corollary}
\begin{proof}
  It is sufficient to approximate $(X_0,V_0)\in \Lspace p$ by a sequence $(X^n_0, V^n_0)\in \hat\cX$
  of initial data arising from finite discrete distributions of space and velocities
  in $\DXspace(\R)$ and to apply the previous Lemma.
\end{proof}
\begin{corollary}[Equivalent characterizations]
  \label{cor:uniqueness}
  Let $(X_0,V_0)\in \Lspace2$ be given initial data. If $X$ is a solution of one
  of the equations \eqref{eq:44}, \eqref{eq:49}, \eqref{eq:26}, then
  it satisfies all the formulations \emph{(L.I, II, III)} and the properties \emph{(L.a,b,c)}
  stated in Theorem \ref{thm:L2flow}.
\end{corollary}
\begin{proof}
  The thesis is obvious in the case of \eqref{eq:44} and \eqref{eq:49}, whose solution is unique
  and it should coincide with the Lagrangian evolution provided by Corollary
  \ref{cor:existence}.

  Let us now assume that $X$ is a Lipschitz curve solving \eqref{eq:26}, let
  $\tilde X$ be the Lagrangian solution given by the previous Corollary \eqref{cor:existence}
  with initial data $X_0,V_0$, 
  and let us set
  $V^n_0:=n(X(n^{-1})-X_0)$, $X^n(t):=\Proj\K(X_0+tV^n_0)$. $X^n(t)$ is thus a Lagrangian flow
  satisfying (L.I, II, III) with respect to the initial data $X_0, V^n_0$; in particular
  \begin{equation}
    \label{eq:70}
    t\frac \d{\dt} X^n(t)\in -\partial I_\K(X^n(t))+X^n(t)-X_0,\quad
    X^n(n^{-1})=X(n^{-1}),
  \end{equation}
  so that $X^n(t)=X(t)$ for $t\ge n^{-1}$. On the other hand, the stability Lemma
  \ref{le:Lagrangian_stability} yields
  \begin{equation}
    \label{eq:105}
    \|X^n(t)-\tilde X(t)\|\le t\|V^n_0-V_0\|=t\|n(X(n^{-1})-X_0)-V_0\|
    \topref{eq:26}\to 0\quad\text{as }n\uparrow+\infty,
  \end{equation}
  so that $X=\tilde X$.
\end{proof}
\section{The continuous sticky particle system in Eulerian coordinates}
\label{sec:limit}
In this section we conclude the proofs of the various theorems of Section \ref{sec:main}.

\begin{proof}[\textbf{Proof of Proposition \ref{prop:metric}}]
Starting from \eqref{eq:12} it is immediate to check that
$\Dist_p$ is a metric on $\Xspace p(\R)$.
Let us check the equivalence characterization \eqref{eq:14}: assuming first
that $\Dist_p(\mu_n,\mu)\to0$ we obviously have $W_p(\rho_n,\rho)\to0$;
since $X_n=X_{\rho_n}\to X=X_\rho$ and $v_n(X_n)\to v(X)$ in $L^p(0,1)$ as $n\uparrow+\infty$,
for a continuous and bounded test function $\zeta:\R\to\R$ we easily get
\begin{equation}
\begin{aligned}
  \lim_{n\uparrow+\infty}\int_\R \zeta(x)v_n(x)\,\d\rho_n(x)&=
  \lim_{n\uparrow+\infty}\int_0^1 \zeta(X_n(w))v_n(X_n(w))\,\d w\\&=
  \int_0^1 \zeta(X(w))v(X(w))\,\d w=
  \int_\R \zeta(x)v(x)\,\d\rho(x),
\end{aligned}\label{eq:38}
\end{equation}
showing that $\rho_nv_n\weakto \rho v$, and
\begin{displaymath}
  \lim_{n\uparrow+\infty}\int_\R |v_n(x)|^p\,\d\rho_n(x)=
  \lim_{n\uparrow+\infty}\int_0^1 |v_n(X_n(w))|^p\,\d w=
  \int_0^1 |v(X(w))|^p\,\d w=\int_\R |v(x)|^p\,\d\rho(x).
\end{displaymath}
The converse implication is a particular case
of \cite[Theorem 5.4.4]{Ambrosio-Gigli-Savare08}: here is a simplified argument.
If \eqref{eq:14} holds, then one gets the strong convergence of $X_n$
to $X$ in $L^p(0,1)$; since $V_n:=v_n\circ X_n$ is bounded in $L^p(0,1)$, up to the extraction
of a suitable subsequence, one has $V_n\weakto V$ in $L^p(0,1)$ and arguing as in \eqref{eq:38}
\begin{equation}
  \label{eq:80}
  \int_0^1 \zeta(X(w))V(w)\,\d w=
  \int_0^1 \zeta(X(w))v(X(w))\,\d w\quad
  \forall\,\zeta\in C_b(\R).
\end{equation}
Notice that a function in $L^p(0,1)$ of the form $b\circ X$ for some Borel map
$b:\R\to \R$ belongs to $\cH_X$; a simple approximation argument shows that
the set $\{\zeta\circ X:\zeta\in C_b(\R)\}$ is dense in $\cH_X$ so that \eqref{eq:80}
yields 
\begin{equation}
  \label{eq:39}
  v\circ X =\Proj{\cH_X}V.
\end{equation}
On the other hand, the last limit property stated in \eqref{eq:14} yields
\begin{equation}
  \label{eq:81}
  \|V\|_{L^p(0,1)}\le \lim_{n\uparrow+\infty}\|V_n\|_{L^p(0,1)}=\|v\circ X\|_{L^p(0,1)}=
  \|\Proj{\cH_X}(V)\|_{L^p(0,1)}\le \|V\|_{L^p(0,1)},
\end{equation}
so that $v\circ X$ should coincide with $V$ which is also the strong limit of $V_n$ in $L^p(0,1)$.

Let us finally consider the density of $\DXspace$: if $(\rho,\rho v)\in \Xspace p(\R)$ we can first approximate
$v$ in $L^p_\rho(\R)$ by a sequence of bounded and continuous functions $v_n\in C_b(\R)$. We can
then find a sequence $\rho^N=\sum_{j=1}^N m_{j,N}\delta_{x_{j,N}}$, $N\in \N$, such that
$\rho^N\to \rho$ in $\cP_p(\R)$. It is then easy to check that $v_n\rho^N\weakto v_n\rho$ as
$N\uparrow+\infty$ according to \eqref{eq:14}.
\end{proof}

\begin{proof}[\textbf{Proof of Theorem \ref{thm:precise_sol}}]
\ \\
\emph{(a)}
The extension of the semigroup $\SS$ is not difficult, by using the estimates of
Theorem \ref{thm:main1}
and the density of $\DXspace(\R)$ in $\Xspace p(\R)$, but not completely trivial since the space
$\Xspace p(\R)$ is not complete and \eqref{eq:17}/\eqref{eq:78} do not provide a
pointwise continuous dependence of the velocity from the initial data. Therefore, we will use
the equivalence stated in Theorem \ref{thm:L2flow}, which we already proved
at the level of discrete data in Theorem \ref{thm:discrete_Lagrangian}, and
the Lagrangian stability result of Lemma \ref{le:Lagrangian_stability}.
It is clear that the only possible extension of $\SS_t$ to $\Xspace p(\R)$
is given by formula \eqref{eq:114}. Since $\sfS_t$ is a semigroup in
$\Lspace p$ satisfying $\lim_{t\downarrow 0}\sfS_t(X_0,V_0)=
(X_0,V_0)$ strongly in $L^p(0,1)^2$, $\SS_t$ satisfies \eqref{eq:110}.

In order to check that $\SS_t$ is strongly-weakly continuous, we
take a sequence
$\mu^n_t=(\rho^n_t,\rho^n_t v^n_t)=\SS_t[\mu^n_0]\in \DXspace$, with
$\mu^n_0$ converging to $\mu=(\rho,\rho v)\in \Xspace p(\R)$
with respect to $\Dist_p$ and we consider the associated monotone rearrangement
maps  $(X^n(t),V^n(t))=\sfS_t(X^n_0,V^n_0)$.
By Lemma \eqref{le:Lagrangian_stability} \emph{(f)}, for every weakly converging sequence $V^{n_k}\weakto \bar V$
in $L^p(0,1)$ and every test function $\zeta\in C^0_b(\R)$ we have
\begin{align*}
  \int_\R \zeta\,v^{n_k}_t\,\d\rho^{n_k}_t&=
  \int_0^1 \zeta(X^{n_k}(t))v^{n_k}_t(X^{n_k}(t))\,\d w
  \topref{eq:74}=
  \int_0^1 \zeta(X^{n_k}_t)V^{n_k}(t)\,\d w
  \\&\stackrel{k\uparrow+\infty}\longrightarrow
  \int_0^1 \zeta(X(t))\bar V\,\d w\stackrel{\text{\it Lemma \ref{le:Lagrangian_stability}(e)}}=
  \int_0^1 \zeta(X(t))V(t)\,\d w\topref{eq:74}=
  \int_\R \zeta\,v_t\,\d\rho_t
\end{align*}
where we used the fact that $\zeta(X^n(t))\to\zeta(X(t))$ strongly in $L^p(0,1)$.

\emph{(b)}
It is immediate to check that  $(\rho,\rho v)=\SS(\rho_0,\rho_0 v_0)$ 
is a distributional solution of \eqref{eq:Euler}, since in Lagrangian coordinates
the continuity equation reads
\begin{align*}
  \frac \d{\d t}\int_\R \zeta(x)\,\d\rho_t(x)&=
  \frac{\d}{\d t}\int_0^1 \zeta(X(t))\,\d w\topref{eq:74}=
  \int_0^1 \zeta'(X(t))V(t)\,\d w\\&\topref{eq:74}=
  \int_0^1 \zeta'(X(t))v_t(X(t))\,\d w=
  \int_\R \zeta'(x)v_t(x)\,\d\rho_t(x),
\end{align*}
and the momentum equation becomes similarly
\begin{align*}
   \frac \d{\d t}\int_\R \zeta(x)v_t(x)\,\d\rho_t(x)&=
   \frac{\d}{\d t}\int_0^1 \zeta(X(t))V(t)\,\d w\stackrel{\eqref{eq:116}\,\eqref{eq:74}}=
   \frac{\d}{\d t}\int_0^1 \zeta(X(t))V_0\,\d w\\&\topref{eq:74}=
   \int_0^1 \zeta'(X(t))V(t)V_0\,\d w=
   \int_0^1 \zeta'(X(t))v_t(X(t))V_0\,\d w\\&\stackrel{\eqref{eq:116}\,\eqref{eq:74}}=
   \int_0^1 \zeta'(X(t))v_t^2(X(t))\,\d w=
   \int_0^1 \zeta'(x)v_t^2(x)\,\d\rho_t(x).
\end{align*}
Oleinik entropy condition \eqref{eq:55} follows easily by \eqref{eq:85}, by observing that
$\Proj{\cH_{X(t)}}(X_0)$ is a nonincreasing map, $V(t)=v_t(X(t))$, and $\rho_t=(X(t))_\#\lambda$.

\emph{(c)} follows from \eqref{eq:106}.

\emph{(d)} is equivalent to point \emph{(d)} of Lemma \ref{le:Lagrangian_stability}; 
concerning the left continuity of $\rho_tv_t$ in the weak topology, we fix an arbitrary bounded 
Lipschitz test function $\zeta:\R\to \R$ and we observe that
\begin{displaymath}
  \lim_{s\uparrow t}\int_\R \zeta(x)v_s(x)\,\d \rho_s(x)=
   \lim_{s\uparrow t}\int_0^1 \zeta(X(s))V(s)\,\d w=
    \lim_{s\uparrow t}\int_0^1 \zeta(X(t))V(s)\,\d w
\end{displaymath}
since $X(s)\to X(t)$ in $L^2(0,1)$ as $s\uparrow t$. On the other hand, since
$\zeta\circ X(t)\in \cH_{X(t)}$ we have
\begin{displaymath}
  \int_0^1 \zeta(X(t))V(s)\,\d w=
  \int_0^1 \zeta(X(t))V(t)\,\d w=
  \int_0^1 \zeta(X(t))v_t(X(t))\,\d w=
  \int_\R \zeta(x)v_t(x)\,\d \rho_t(x).
\end{displaymath}

\emph{(e)} has already been discussed in point \emph{(a)}, except for the convergence
at $t\in (0,+\infty)\setminus \cT$, which follows from Lemma \ref{le:Lagrangian_stability} \emph{(f)}.

\emph{(f)}
\eqref{eq:71} follows by the projection representation \eqref{eq:83}
and Corollary \ref{cor:omega}. The limit in \eqref{eq:71} can be obtained in Lagrangian coordinate:
\begin{align*}
  \lim_{h\downarrow 0}\int_\R \Big|h^{-1}(\sft_s^{s+h}-\sfi)-v_s\Big|^2\,\d\rho_s=
  \lim_{h\downarrow 0}\int_0^1 \Big|h^{-1}(X(s+h)-X(s))-V(s)\Big|^2\,\d w=0
\end{align*}
since $t\mapsto X(t)$ is right differentiable.
\eqref{eq:72} is an immediate consequence of \eqref{eq:85}, which yields
\begin{displaymath}
  (t-s)V(t)=X(t)-\Proj{\cH_{X(t)}}(X(s))\quad
  \forall\, 0\le s<t.\qedhere
\end{displaymath}
\end{proof}
\begin{proof}[The \textbf{proof of Theorem \ref{thm:L2flow}}]
follows now by applying Lemma \ref{le:Lagrangian_stability} and
its corollaries \ref{cor:existence}, \ref{cor:uniqueness}.
\end{proof}
\begin{proof}[\textbf{Proof of Theorem \ref{thm:WGFlow}}]
\eqref{eq:24}
follows from a simple calculation starting from \eqref{eq:26}: we 
introduce the monotone rearrangement $Z$ of the measure $\eta\in \cP_2(\R)$ and we 
observe that $W_2^2(\rho_t,\eta)=\|X(t)-Z\|^2$ (we use the usual notation for $(X,V)$
and we denote by $\|\cdot\|$ the norm
in $L^2(0,1)$).
We get for some $\Xi(t)\in \partial I_\K(X(t))$
\begin{align*}
  \label{eq:41}
  \frac t2\Rightdert W_2^2(\rho_t,\eta)&=
  \frac t2 \Rightdert \|\X(t)-Z\|^2=
  t\Scalar{ \dot \X(t)}{\X(t)-Z}
  \topref{eq:26}=
  \Scalar{\X(t)-\X_0-\Xi(t)}{\X(t)-Z}\\&\topref{eq:117}\le
  \Scalar{\X(t)-\X_0}{\X(t)-Z}
  =
  \frac 12  \|\X(t)-Z\|^2 -\frac 12 \|Z-\X_0\|^2
  +\frac 12 \|\X(t)-\X_0\|^2
  \\&=\frac 12 W_2^2(\rho_t,\eta)-\phi^{\rho_0}(\rho_t)+\phi^{\rho_0}(\eta).
\end{align*}
Let us consider now the converse implication: if $\rho_t$ satisfies \eqref{eq:22} then
$X(t)=X_{\rho_t}$ satisfies (see \eqref{eq:22})
\begin{equation}
  \label{eq:118}
  \frac t2\frac \d{\dt}\|X(t)-Z\|^2-\frac 12\|X(t)-Z\|^2\le \Phi^{\rho_0}(Z)-\Phi^{\rho_0}(X(t))
  \quad\forall\, Z\in \K,
\end{equation}
which is the equivalent metric formulation \cite{Ambrosio-Gigli-Savare08} of the differential inclusion \eqref{eq:26}.

Since $\rrho_t=(X_0,X(t))_\#\lambda$, \eqref{eq:109} yields
\begin{equation}
  \label{eq:119}
  \lim_{t\downarrow0}t^{-2}\int_0^1 |X_0+tV_0-X(t)|^2\,\d w=0,
\end{equation}
i.e. $X(t)$ also satisfies the initial limit condition of \eqref{eq:26}. Therefore, setting $V:=\frac \d\dt X=v\circ X$,
by Corollary
\ref{cor:uniqueness} the couple $(X(t), V(t))$ coincides with the Lagrangian flow $\sfS_t(X_0,V_0)$
so that $(\rho_t,\rho_t v_t)=\SS_t(\rho_0,\rho_0 v_0).$ 
\end{proof}
\begin{proof}[\textbf{Proof of Theorem \ref{thm:limit}}]
  Let us first notice that when $\sfi+\eps_0v_0$ is $\rho_0$-essentially nondecreasing,
  \eqref{eq:19}
  follows directly from \eqref{eq:28}, since the collision-free motion
  $\rho_t=(\sfi+t v_0)_\#\rho_0$ for $t\in [0,\eps_0)$ is a solution of the sticky particle system.

  Let us now consider the general case, setting
  $\tilde\rho_{\eps, t}:=\GG^{\rho_0}_{\log (t/\eps)}(\tilde\rho_\eps)$.
  For every $\eps>0$ let us consider the convex set of bounded Lipschitz functions
  \begin{displaymath}
    BL(\eps):=\Big\{u\in C^{0,1}(\R): \sup|u|\le \eps^{-1}, {\rm Lip}(u)\le (2\eps)^{-1}\Big\}      
  \end{displaymath}
  and let $u_\eps\in BL(\eps)$ be a minimizer of
  \begin{equation}
    \label{eq:120}
    m_\eps=\min_{u\in BL(\eps)}\|v_0-u\|=\|v_0-u_\eps\|.
  \end{equation}
  By standard approximation results, $\lim_{\eps\downarrow0} m_\eps=0$,
  so that $u_\eps$ converges to $v_0$.

  By the definition of $BL(\eps)$ the map $\sfi+\eps u_\eps$ is monotone, and therefore it is the optimal
  map pushing $\rho$ to $\hat \rho_\eps=(\sfi+\eps u_\eps)_\#\rho_0$.
  The sticky particle solution $(\hat\rho_{\eps,t},\hat \rho_{\eps, t}\hat v_{\eps,t}):=\SS_t(\hat \rho_0,\hat \rho_0 u_\eps)$
  admits the representation \eqref{eq:28}
  \begin{displaymath}
    \hat\rho_{\eps,t}=\GG^{\rho_0}_{\log(t/\eps)}(\hat \rho_\eps)
  \end{displaymath}
  so that, by the exponential rate of expansion of $\GG$ we get
  \begin{equation}
    \label{eq:121}
    W_2(\hat\rho_{\eps,t},\tilde\rho_{\eps,t})\le \exp\big(\log (t/\eps)\big) W_2(\hat\rho_\eps,\tilde\rho_\eps)=
    \frac t\eps W_2(\hat\rho_\eps,\tilde\rho_\eps)\le t\|v_0-u_\eps\|_{L^2_{\rho_0}(\R)}\topref{eq:120}=t m_\eps.
  \end{equation}
  On the other hand, if $(\rho_t,\rho_t v_t)=\SS_t(\rho_0,\rho_0v_0)$, \eqref{eq:16} yields
  \begin{equation}
    \label{eq:122}
    W_2(\hat\rho_{\eps,t},\rho_t)\le t\|v_0-u_\eps\|_{L^2_{\rho_0}(\R)}=t m_\eps,
    \quad\text{so that}\quad
    W_2(\rho_{t},\tilde\rho_{\eps,t})\le 2m_\eps t,
  \end{equation}
  and concludes the proof of \eqref{eq:29}.
\end{proof}

We conclude this section by showing that the representation-convergence Theorem
of \textsc{Brenier \band Grenier} \cite{Brenier-Grenier98} can be easily deduced
by our result, in particular by formula \eqref{eq:49} of Theorem \ref{thm:L2flow}.
\begin{theorem}[Brenier-Grenier]
  \label{thm:Brenier-Grenier}
  Let $v_0\in C^0(\R)$, let $\rho_0^N$, $N\in \N$, be a sequence
  of discrete probability measures supported in a fixed compact interval $[-R,R]$ and
  weakly converging to $\rho_0$ in $\cP(\R)$, and let
  $\rho^N_t$ be the solution of the discrete SPS with initial
  data $(\rho_0^N, v_0\rho_0^N)$.
  For every $t\ge0$ $\rho^N_t$ weakly converge to a probability measure $\rho_t$,
  whose distribution function $M_t(x):=\rho_t((-\infty,x])$, $t\ge 0$,
  is the unique entropy solution of
  \begin{equation} \label{BG}
    \partial_t M + \partial_x(A(M)) = 0, \quad M(0)= M_0,
  \end{equation}
  where the flux function $A:[0,1] \to\R$ is defined by
  \begin{equation}
    \label{eq:123}
    A(w):=\int_0^w V_0(r)\,\d r,\quad \text{where}\quad
    V_0:=v_0\circ X_0,\quad X_0:=X_{\rho_0}.
  \end{equation}
\end{theorem}
\begin{proof}
  The convergence part follows by Theorem \ref{thm:precise_sol} and we can represent 
  $X_t:=X_{\rho_t}$ by the formula $X_t=\Proj\K(X_0+tV_0)$ of Theorem \ref{thm:L2flow}.
  Introducing the
  convex primitive functions $F_t(w):=\int_0^w X_t(r)\,\d r$,
  Theorem \ref{proj on K} yields
  \begin{equation}
    \label{eq:125}
    F_t=\big(F_0+tA\big)^{**}\quad\text{so that}\quad \big(F_t\big)^*=\big(F_0+tA\big)^*.
  \end{equation}
  On the other hand, since the derivative $X_t$ of $F_t$ is the pseudoinverse of $M_t$ \eqref{eq:90},
  a standard duality result shows that
  $\big(F_t\big)^*=G_t$ where $G_t(x)=\int_{-\infty}^x M_t(y)\,\d y$, so that
 \begin{equation}
    \label{eq:126}
    G_t=\big(F_0+tA\big)^*=\big(G_0^*+tA\big)^*
  \end{equation}
  It was already observed by \cite[\S 4]{Brenier-Grenier98} that \eqref{eq:126} provides
  the second Hopf formula \cite{Bardi-Evans84} for the viscosity solution of the Hamilton-Jacobi equation
  \begin{equation}
    \label{eq:127}
    \partial_t G+A(\partial_x G)=0\quad\text{in }\R\times (0,+\infty),
  \end{equation}
  and therefore the derivative $M_t=\partial_x G_t$ is the entropy solution of \eqref{BG}.
\end{proof}

\def\cprime{$'$} \def\cprime{$'$}

\end{document}